\documentclass[10pt,reqno]{amsart}
\usepackage{amssymb,mathrsfs,color}
\allowdisplaybreaks
\usepackage{enumerate}
\usepackage{graphicx}
\usepackage{xcolor} 
\usepackage{geometry}
\usepackage{amsfonts,amssymb,amsmath,mathtools}
\usepackage{slashed}
\usepackage{ mathrsfs }
\usepackage[titletoc,title]{appendix}
\usepackage{wrapfig}

\usepackage{cancel}

\usepackage{cite}

\definecolor{green}{rgb}{0,0.8,0} 

\title[Large time behavior]{Large time behavior of global solutions to one-dimensional non-isentropic compressible Navier-Stokes-Poisson equations with large initial data}

\author[Zhang]{Lan Zhang}
\address[Lan Zhang]{\newline School of Mathematics and Statistics and Center for Mathematical Sciences,
Wuhan University of Technology, Wuhan 430070, China\\
\newline Hubei Key Laboratory of Mathematics \& Intelligent Computation and Interdisciplinary Sciences, Wuhan University of Technology, Wuhan 430070, China}
\email{zhanglan91@whut.edu.cn}

\author[Zhao]{Huijiang Zhao}
\address[Huijiang Zhao]{\newline School of Mathematics and Statistics, Wuhan University, Wuhan 430072, China\\
 \newline Computational Science Hubei Key Laboratory, Wuhan University, Wuhan 430072, China}
\email{hhjjzhao@hotmail.com}

\newcommand{\bbR}{\mathbb R}

\definecolor{deepgreen}{cmyk}{1,0,1,0.5}

\makeatletter

\newcommand{\Rmnum}[1]{\expandafter\@slowromancap\romannumeral #1@}
\makeatother

\newcommand{\Del}[1]{}

\numberwithin{equation}{section}

\newtheorem{theorem}{Theorem}[section]
\newtheorem{corollary}[theorem]{Corollary}
\newtheorem{lemma}[theorem]{Lemma}
\newtheorem{remark}[theorem]{Remark}

\begin{document}
\date{\today}

\thanks{
\textbf{Acknowledgment.} The research of Lan Zhang was supported by the grants from National Natural Science Foundation of China (Grants No. 12101472 and 12471210).
The research of Huijiang Zhao is supported by the National Natural Science Foundation of China (Grants No. 12221001, 12371225, and 12571242).}


\begin{abstract}
In this paper, we are concerned with the construction of global solutions to the Cauchy problem for the one-dimensional non-isentropic compressible Navier-Stokes-Poisson equations, as well as the precise description of their large-time behavior. While results for small initial data are well-established, the case of large initial perturbations remains less satisfactory. Even for constant transport coefficients, only the Nishida-Smoller type solutions have been obtained so far. In this work, we will provide a complete resolution of this problem. The key idea involves revealing favorable hidden identities of the potential force and a new representation for the specific volume, which enables us to derive uniform upper and lower bounds on this quantity. Moreover, this strategy extends to compressible viscous single-ion flow under the Boltzmann relation; notably, this new method for deriving lower and upper bounds on the specific volume also applies to the Navier-Stokes equations.
The proof relies on a domain decomposition technique and some new dissipative structures to close the nonlinear energy estimates.
\\[2mm]

\\[2mm]
\end{abstract}

\maketitle \centerline{\date}



\section{Introduction}
\subsection{The models}
The compressible Navier-Stokes-Poisson equations (denote as NSP in the sequel) which take form of compressible Navier-Stokes equations
coupled with self-consistent Poisson equation are often used to simulate the motion of viscous fluid under the influence
of self-consistent electrostatic potential force in semiconductor devices
or plasma physics, c.f. \cite{chen1984introduction,MR1063852}.
In this paper, we will consider the one-dimensional non-isentropic compressible NSP equations in Lagrangian coordinates as follows:
\begin{equation}\label{nsp}
  \begin{cases}
  v_t-u_{x}=0,\\
  u_t+p(v,\theta)_x=\left(\frac{\mu u_x}{v}\right)_x+\frac{\Phi_x}{v},\\
  e_t+p(v,\theta)u_x=\frac{\mu u^2_x}{v}+\left(\frac{\kappa \theta_x}{v}\right)_x,\\
  \left(\frac{\Phi_x}{v}\right)_x=1-ve^{-\chi\Phi}-vb(x),\quad \lim_{|x|\rightarrow+\infty}\Phi(t,x)=0.
  \end{cases}
\end{equation}
Here, the unknown functions $v(t,x)$, $u(t,x)$ and $\theta(t,x)$ stand for the specific volume, velocity and temperature,
respectively. $\Phi(t,x)$ is the self-consistent potential for any $(t,x)\in\bbR_+\times\bbR$. The viscous coefficient $\mu$ and heat-conductivity coefficient $\kappa$ are assumed to be
generic positive constants or to depend on $v(t,x)$ and $\theta(t,x)$.
Moreover, we will consider the ideal, polytropic gases:
\begin{equation}\label{pressure}
  p(v,\theta)=\frac{R\theta}{v}=Av^{-\gamma}\exp\left(\frac{\gamma-1}{R}s\right),\quad e=C_v\theta=\frac{R\theta}{\gamma-1}.
\end{equation}
Here $A$ is a positive constant, $s$ is entropy, the positive constants $R$ and $C_v$ are specific gas constant and specific heat constant, respectively. The adiabatic constant $\gamma>1$.

The given function $b(x)$ is the doping profile, which stands for the
density of impurities in charged fluid or semiconductor devices. In this paper, $b(x)$ is supposed to be $0$.
The given parameter $\chi$ is artificially assumed, and it takes the values $0$ and $1$ depending on the physical assumptions made about the electron dynamics.
In the two-fluid model for charged particles where the Poisson equation satisfies $\Phi_{xx}=n_e-n_i$ with $n_e$ and $n_i$ being the density of electrons and ions (cf. \cite{MR3306347,MR4458894,MR4362505}), thanks to the fact that the mass of an electron is much smaller than that of an ion, the much heavier ions can be treated as motionless, with a constant density $n_i = n_0$ (which is normalized to $1$), forming only a fixed charged background. Then the electrostatic potential $\Phi$ is determined by this framework, leading to the Poisson equation formulated as $\left(\frac{\Phi_x}{v}\right)_x=1-v$, $v$ is the specific volume of the electrons. Such a system is used for describing Langmuir waves (electron waves) in the two-fluid theory (cf. \cite{MR1671928}).
On the other hand, under the assumption that the electrons are massless,
they respond instantaneously to the electric field and reach thermodynamic equilibrium in much faster rate than heavier ions in plasma (cf. \cite{chen1984introduction,MR2825596,MR2775116}).
This gives rise to the celebrated Boltzmann relation $n_e=e^{-\Phi}$ for the electron density,
which yields $\left(\frac{\Phi_x}{v}\right)_x=1-ve^{-\Phi}$, $v$ is the specific volume of the ions.
This model, often referred to as the one-fluid model for ions, captures the essential nonlinear physics of electrostatic shielding and is particularly effective for studying large-scale, low-frequency ion dynamics where the rapid equilibration of electrons can be reasonably assumed.

It is well-known that the compressible Euler-Poisson equations admit global-in-time smooth solutions with small amplitude (cf. \cite{MR4362505,MR1671928,MR2775116,MR3595365}). However, thanks to the quasi-linear hyperbolic nature, as expected, shock waves do develop for {\it large} perturbations of the constant state equilibrium (cf. \cite{MR1724661}). To suppress the shock formation for large initial perturbation, the known results considered the viscous effects,
and through these effects the global existence and large time behavior of solutions to one-dimensional compressible Navier-Stokes equations have been derived (cf. \cite{MR0227619,MR0468593,MR1671920,MR1912419,MR3466845}).
We wonder whether it is valid for one-dimensional compressible Navier-Stokes-Poisson equations.
Recently, the authors \cite{MR4458894} have proven such type results for the isentropic case. However, the method used there to
obtain uniform bounds for the specific volume is specific to the isentropic case and fails entirely
in the non-isentropic setting, where the evolution of temperature and its coupling with the energy
equation introduce substantial new difficulties. In this paper, we will focus on this open question and develop some new ideas.

\subsection{The main results}
The main purpose of this paper is concerned with the construction of global strong solutions and the precise description
of the large time behaviors of the global solutions constructed to the Cauchy problem of the system \eqref{nsp} with
the following initial data:
\begin{equation}\label{initial-data}
  \begin{aligned}
  (v,u,\theta)(0,x)&=(v_0,u_0,\theta_0)(x),\\
  \lim_{x\rightarrow\pm\infty}(v_0,u_0,\theta_0)(x)&=(v_{\pm},u_{\pm},\theta_{\pm}).
  \end{aligned}
\end{equation}
Without loss of generality, we assume $v_-=v_+=1$, $u_-=u_+=0$, $\theta_-=\theta_+=1$.

\medskip
In our first theorem, we will consider the Cauchy Problem \eqref{nsp}-\eqref{initial-data} with the transport coefficients being generic positive constants. The main results are as follows:
\begin{theorem}\label{main-theorem}
Assume that the initial data $(v_0,u_0,\theta_0)$ defined in \eqref{initial-data} satisfy
\begin{align}\label{initial-assump}
  (v_0-1,u_0,\theta_0-1)\in H^1(\bbR),\quad
  0<\underline{V}\leq v_0(x)\leq \overline{V},\quad \inf_{x\in\bbR}\theta_0(x)>0,
\end{align}
where $\overline{V}>1$ is a given positive constant.
When $\chi=0$, we assume further that $\Phi_{0x}\in H^1(\bbR)$.
Then there exists a unique global solution $(v,u,\theta,\Phi)(t,x)$ to system \eqref{nsp}-\eqref{initial-data} for any $T>0$ such that
\begin{itemize}
  \item When $\chi=0$,
  \begin{equation}\label{regularity}
  \begin{cases}
  (v-1,u,\theta-1)\in L^{\infty}(0,T;H^1(\bbR)),\quad \Phi_{x}\in L^{\infty}(0,T;H^1(\bbR));\\
  (v_x,1-v)\in L^2(0,T;L^2(\bbR)),\quad (u_x,\theta_x)\in L^2(0,T;H^1(\bbR)).
  \end{cases}
\end{equation}
  \item When $\chi=1$,
  \begin{equation}\label{regularity-2}
  \begin{cases}
  (v-1,u,\theta-1)\in L^{\infty}(0,T;H^1(\bbR)),\quad \Phi\in L^{\infty}(0,T;H^2(\bbR))\\
  v_x\in L^2(0,T;L^2(\bbR)),\quad (u_x,\theta_x,\Phi_x)\in L^2(0,T;H^1(\bbR)).
  \end{cases}
\end{equation}
\end{itemize}
Moreover, there exists constants $C_{i}(i=0\sim6)$ which depend only on $\mu$, $\kappa$, $R$, $C_v$ and initial assumptions \eqref{initial-assump}
such that
\begin{itemize}
  \item When $\chi=0$,
  \begin{gather}
  C_0^{-1} \leq v(t,x) \leq C_0, \quad
  C_1^{-1} \leq \theta(t,x) \leq C_1, \\
  \sup_{0\leq t\leq T}\left\{\|(v-1,u,\theta-1)(t)\|_{1}^2+\|\Phi_x(t)\|_{1}^2\right\}
    +\int_0^T\|(1-v,u_x,\theta_{x})(t)\|_1^2dt
    \leq C_2,\\
  \lim_{t\rightarrow\infty}\left(\|(v-1,u,\theta-1,\Phi_x)\|_{L^p}+\|(1-v,v_x,u_x,\theta_x)\|_{L^2}\right)=0,\quad {\rm for}\ \forall p\in(2,+\infty].
\end{gather}
  \item When $\chi=1$,
  \begin{gather}
  C_3^{-1} \leq v(t,x) \leq C_3, \quad
  C_4^{-1} \leq \theta(t,x) \leq C_4, \quad
  |\Phi(t,x)|\leq C_5,\\
  \sup_{0\leq t\leq T}\left\{\|(v-1,u,\theta-1)(t)\|_{1}^2+\|\Phi(t)\|_{2}^2\right\}
    +\int_0^T\left(\|v_x\|_{L^2}^2+\|(u_x,\theta_{x},\Phi_x)(t)\|_1^2\right)dt
    \leq C_6,\\
  \lim_{t\rightarrow\infty}\left(\|(v-1,u,\theta-1,\Phi,\Phi_x)\|_{L^p}+\left\|\left(1-ve^{-\Phi},v_x,u_x,\theta_x,\Phi_x\right)\right\|_{L^2}\right)=0, \quad {\rm for}\ \forall p\in(2,+\infty].
\end{gather}
\end{itemize}
\end{theorem}

\medskip
Our next theorem is concerned with a physical case which means the viscosity coefficient $\mu(v,\theta)$ and heat-conductivity coefficient $\kappa(v,\theta)$ depend on the specific
volume $v(t,x)$ and the temperature $\theta(t,x)$. Actually, NSP system can be transformed from Vlasov-Poisson-Boltzmann equations by applying Chapman-Enskog expansion, see \cite{MR1148892,MR1307620,MR0156656},
which yields the density and temperature dependent transport coefficients $\mu$ and $\kappa$. Up to date, there are lots of papers concerned about the one-dimensional compressible Navier--Stokes
equations with density and temperature dependent transportation coefficients in various form, see \cite{MR2644363,MR791841,MR653464,MR3461630,MR3291375} and the references therein.
For the system \eqref{nsp}, Tan, et al. \cite{MR3032988} considered a case that
\begin{equation}\label{tan}
  \mu(v)=v^{-a},\quad \kappa(v,\theta)=\theta^b,
\end{equation}
where $a$ and $b$ satisfy complex relations. They required the smallness of the adiabatic constant $\gamma-1>0$ to derive the large time behavior of the global solutions with large initial data.
Hence, we wonder whether this smallness assumption can be removed.

We will assume throughout Section 3 that the viscosity coefficient $\mu$ and heat conductivity
coefficient $\kappa$ are smooth functions of the temperature $\theta$ and the specific volume $v$ as follows,
\begin{equation}\label{transport}
\mu=\tilde{\mu}h(v)\theta^{\alpha},\quad
\kappa=\tilde{\kappa}h(v)\theta^{\alpha},
\end{equation}
where $\tilde{\mu}$ and $\tilde{\kappa}$ are positive constants,
and there exist positive constants $C$, $\ell_1$, and $\ell_2$
such that
\begin{equation}\label{h}
h(v)\thicksim v^{\ell_1}+v^{-\ell_2},
\quad h'(v)^2v\leq C h(v)^3 \quad {\rm for\ all}\ v\in(0,\infty).
\end{equation}
Here, the assumptions \eqref{h} are used to control the highly nonlinear terms such as \eqref{3.4-1} and to derive upper and lower bounds of the specific volume $v$ which depend on the temperature $\theta$.
Moreover, it is clear that if we choose $\ell_1=\ell_2=\alpha=0$ the transport coefficients become positive constants which is a special case of \eqref{transport}.
\begin{theorem}\label{main-theorem-1}
Assume that the viscosity $\mu$ and the heat conductivity $\kappa$ satisfy \eqref{transport}--\eqref{h} for
\begin{equation}\label{ell-assump}
  \frac{2}{1+2\ell_1}+\frac{1}{\ell_2}< \frac{4}{5},\quad \ell_1\geq 1,\quad \ell_2\geq 1.
\end{equation}
Let the initial data $(v_0,u_0,\theta_0,\Phi_{0x})$ defined in \eqref{initial-data} satisfy
\begin{align}\label{initial-assump-1}
  (v_0-1,u_0,\theta_0-1)\in H^3(\bbR),\quad
  0<\overline{V}^{-1}\leq v_0(x)\leq \overline{V},\quad \theta_0(x)>\overline{V}^{-1},
\end{align}
where $\overline{V}>1$ is a given positive constant.  When $\chi=0$, we assume further that $\Phi_{0x}\in H^3(\bbR)$. Then there exists a positive constant $\epsilon_0$ such that $|\alpha|\leq \epsilon_0$, the Cauchy problem \eqref{nsp}-\eqref{initial-data}
admits a unique solution  $(v, u, \theta,\Phi)(t,x)$ satisfying
\begin{itemize}
  \item When $\chi=0$,
  \begin{equation}\label{regularity-1}
  \begin{cases}
  (v-1,u,\theta-1)\in L^{\infty}(0,T;H^3(\bbR)),\quad \Phi_{x}\in L^{\infty}(0,T;H^3(\bbR)));\\
  v_x\in L^2(0,T;H^2(\bbR)),\quad (u_x,\theta_x)\in L^2(0,T;H^3(\bbR)).
  \end{cases}
\end{equation}
  \item When $\chi=1$,
  \begin{equation}\label{regularity-3}
  \begin{cases}
  (v-1,u,\theta-1)\in L^{\infty}(0,T;H^3(\bbR)),\quad \Phi\in L^{\infty}(0,T;H^4(\bbR)));\\
  v_x\in L^2(0,T;H^2(\bbR)),\quad (u_x,\theta_x,\Phi_x)\in L^2(0,T;H^3(\bbR)).
  \end{cases}
\end{equation}
\end{itemize}
Moreover, there exists constants $B_{i}(i=0\sim6)$ which depend only on $\mu$, $\kappa$, $R$, $C_v$ and initial assumption \eqref{initial-assump-1} such that
\begin{itemize}
  \item When $\chi=0$,
  \begin{gather}
  B_0^{-1} \leq v(t,x) \leq B_0, \quad
  B_1^{-1}\leq \theta(t,x) \leq B_1, \\
  \sup_{0\leq t\leq T}\left\{\|(v-1,u,\theta-1)(t)\|_{3}^2+\|\Phi_x(t)\|_{3}^2\right\}
    +\int_0^T\|(1-v,u_x,\theta_{x})(t)\|_3^2dt
    \leq B_2,\\
  \lim_{t\rightarrow\infty}\left(\|(v-1,u,\theta-1,\Phi_x)\|_{L^p}+\|(1-v,v_x,u_x,\theta_x)\|_{L^2}\right)=0, \quad {\rm for}\ \forall p\in(2,+\infty].
\end{gather}
  \item When $\chi=1$,
  \begin{gather}
  B_3^{-1} \leq v(t,x) \leq B_3, \quad
  B_4^{-1}\leq \theta(t,x) \leq B_4, \quad
  |\Phi(t,x)|\leq B_5,\\
  \sup_{0\leq t\leq T}\left\{\|(v-1,u,\theta-1)(t)\|_{3}^2+\|\Phi(t)\|_{4}^2\right\}
    +\int_0^T\left(\|v_x\|^2_{2}+\|(u_x,\theta_{x},\Phi_x)(t)\|_3^2\right)dt
    \leq B_6,\\
  \lim_{t\rightarrow\infty}\left(\|(v-1,u,\theta-1,\Phi,\Phi_x)\|_{L^p}+\|(1-ve^{-\Phi},v_x,u_x,\theta_x,\Phi_x)\|_{L^2}\right)=0, \quad {\rm for}\ \forall p\in(2,+\infty].
\end{gather}
\end{itemize}
\end{theorem}

\begin{remark}
In the proof of Theorem \ref{main-theorem}, we do not need any additional smallness assumptions which
extend and improve the previous results such as \cite{MR3032988,MR3436996,MR3466845,MR4884275,Zhao-Zhao-2026,MR4693518,MR4755045}.
For the Theorem \ref{main-theorem-1}, we use the smallness of $|\alpha|$ instead of the Nishida-Smoller type assumption which means the adiabatic constant $\gamma-1$ being sufficiently small which improves the results in \cite{MR3032988,MR4884275}. Notice that we have considered a different set of assumptions on the transport coefficients depending on density and temperature for the system \eqref{nsp} in \cite{MR4693518}, however, which are not covered by the assumptions \eqref{transport}-\eqref{h}.
\end{remark}
\begin{remark}
To derive the upper bound for the specific volume, we combine a spatial antiderivative with averaging over $[x-l,x+l]$, where $l>0$ is chosen
independently of $t$ and $x$. For the lower bound, we apply Jiang's local representation \cite{MR1671920,MR1912419}, but on time intervals of a fixed length $h>0$.
The resulting lower bound depends only on the initial data and $h$, and is independent of the total existence time.
Together with the a priori estimates, this enables us to extend the local solution to a global one. It is worth noting that both methods also
apply to the classical Navier--Stokes equations. Indeed, considering the \eqref{v_form2} with
\begin{equation*}
  A_l(t,x):=\exp\left\{-\frac{R}{\mu}\int_0^t\left\langle\frac{\theta}{v}\right\rangle_{l,x}ds\right\},
\end{equation*}
we can take $l=1$ such that
\begin{align*}
  A_1(t,x)&\leq \exp\left\{-\frac{R}{2\mu}\int_0^t\inf_{y\in[x-1,x+1]}\theta(t,y)\left(\int_{x-1}^{x+1}v(t,y)dy\right)^{-1}ds\right\}\\
  &\leq \exp\left\{-\frac{R}{4\alpha_2\mu}\int_0^t\inf_{y\in[x-1,x+1]}\theta(t,y)ds\right\}.
\end{align*}
Repeating the same estimates in \eqref{low-theta-good} and \eqref{A-l}, we can get
\begin{equation*}
  A_1(t,x)\leq Ce^{-Ct},
\end{equation*}
which yields the upper bound of $v(t,x)$. We believe that this method can be also applied to other hydrodynamic models.
\end{remark}
\begin{remark}
Due to the effect of electric potential force, we could obtain the asymptotic behavior for $\|v-1\|_{L^2}$ in the Theorem \ref{main-theorem} and Theorem \ref{main-theorem-1},
which is absent for the compressible Navier-Stokes equations.
Moreover, an extra uniform bound for the self-consistent electric potential $\Phi(t,x)$ when $\chi=1$
is derived through the damping effect induced by the Boltzmann relation $n_e=e^{-\Phi}$. The large time behavior of terms related to the electric potential such as
$\|\Phi(t,x)\|_{L^p}$ $(p>2)$ and $\|\Phi_x(t,x)\|_{L^2}$ are obtained as well.
\end{remark}
\begin{remark}
It can be anticipated that when the doping profile $b(x)\neq0$, the system of equations will exhibit rich wave phenomena. Recall that the stability of steady states
of multidimensional Navier-Stokes-Poisson equations with non-flat doping profile were established in \cite{MR3296606,MR4941832}.
However, for one-dimensional case with large initial perturbations, even if the far field states of $(v,u,\theta)(t,x)$ are the same, the doping profile $b(x)$ causes great difficulties.
It is readily to check that the hidden relations \eqref{new-relation-0}
and \eqref{new-relation} we found are no longer holds which is crucial to obtain the uniformly upper and lower bounds for the specific volume. We will propose this
difficult problems in the future.
\end{remark}
\begin{remark}
We believe that the analysis presented here can shed light on the construction of basic wave patterns for the non-isentropic compressible Navier-Stokes-Poisson system with large initial perturbations when $\chi=1$. It is worth pointing out that the nonlinear stability of composite wave patterns with large initial data for system \eqref{nsp}
with $\chi=1$ has been established in \cite{Zhao-Zhao-2026}, where the transport coefficients are assumed to depend on the specific volume and temperature. See also other results on the stability of rarefaction waves for a two-fluid full compressible Navier-Stokes-Poisson system under large initial perturbations (cf. \cite{MR4884275}), which required the adiabatic constant $\gamma-1$ to be sufficiently small. For the case of constant transport coefficients, the problem without smallness assumption remains open.
\end{remark}

\subsection{Former results related and our main ideas}
There are many results concerning the large time behavior of global solutions to conservation laws with large initial data.
Especially, for one dimensional compressible Navier--Stokes equations, Kanel$'$ \cite{MR0227619}
first studied the isentropic case with $p=Rv^{-\gamma}$ for $\gamma>0$ and obtained the existence and large time behavior of global solutions for large initial data.
For the non--isentropic case, the existence of global solutions in bounded domain with large initial data was first obtained by Kazhikhov and Shelukhin in \cite{MR0468593}. After that, many efforts have been
made to show the asymptotic behavior of global solution as $t\rightarrow\infty$. With some additional smallness conditions on the initial data, the corresponding results are confirmed
in \cite{MR1142276,MR1748226,MR530493,MR637519,MR1357824,MR692729} and the references therein. To remove the smallness assumptions, Jiang \cite{MR1671920,MR1912419} firstly
applied a local representation of the specific volume to prove its pointwise upper and lower bounds independent of time and space, with which
the global solutions are convergent locally in space as time goes to infinity for large initial data and the temperature is bounded from below and above locally in space for all $t\geq0$.
Whether the global solution converges to a constant steady state uniformly with respect to the spatial variable as time goes to infinity for large initial data, Li-Liang \cite{MR3466845} gave a positive answer
by deriving lower and upper bounds for the temperature independent of both time and space. Then, Wang-Zhao \cite{MR3564590} extended the corresponding results to the case that transport coefficients depend on
the density and temperature and are both proportional to $h(v)\theta^\alpha$ for certain non-degenerate smooth function $h(v)$.
From then on, there are a lot of recent papers on the stability of wave patterns to the one-dimensional compressible Navier--Stokes equations with large initial perturbation, see \cite{MR3861291,MR3703560,MR3595482,MR3794193} and the references therein.

When the effect of potential force is taken into account, there are few results on the existence and asymptotic behavior of global solutions to the compressible NSP equations with large initial data. The inclusion of the Poisson potential introduces a nonlocal term $\frac{\Phi_x}{v}$
that fundamentally changes the mathematical structure. This term prevents the direct application of existing techniques in \cite{MR1671920,MR0468593,MR3466845}
 to obtain uniform-in-time bounds on the specific volume $v$
and temperature $\theta$, especially for large initial data. As a result, the large-time behavior of solutions to the non-isentropic compressible NSP system with general large data has remained a challenging problem.

For one-dimensional non-isentropic NSP system, Tan-Yang-Zhao-Zou \cite{MR3032988} studied the global solutions with large initial data away from vacuum to the Cauchy problem
with density-dependent viscosity coefficient and density- and temperature-dependent heat-conductivity coefficient. Moreover, they obtained the large time behavior of the global solutions in the sense of Nishida-Smoller type (cf. \cite{MR330789}) when
the adiabatic constant $\gamma-1$ is sufficiently small. Recently, \cite{MR4693518} removed the smallness assumption by assuming transport coefficients depending on
both density and temperature but with complex constraint conditions. Regarding the isentropic case, we refer to the recent work  \cite{MR4755045} by the first author.
Considering the nonlinear stability of wave patterns,
the authors \cite{MR4458894} first showed that rarefaction waves of the compressible Navier-Stokes-Poisson system are time-asymptotically nonlinear stable both for the two-fluid case of charged particles consisting of ions and electrons and for the one-fluid case of single ions under the Boltzmann relation. After that, the asymptotic
stability of the shock wave, rarefaction wave and composite wave as well as nonlinear wave for inflow/outflow problem were accomplished in \cite{MR4903150,MR4547899,MR4610411,MR4842676,Zhao-Zhao-2026,MR4884275} and the references therein.
For higher dimensional case, Liu-Luo-Zhong \cite{MR4122651} considered the initial boundary value problems for compressible NSP system on the exterior domain. With the radial symmetry assumption,
the global existence of solutions with large initial data on a domain exterior to a ball in $\bbR^n(n\geq1)$ is proved.
Concerning an initial value problem of 3D compressible NSP equations subject to large and non-flat doping profile,
\cite{MR4113208,MR4444174} established the global well-posedness of strong solutions with large oscillations and vacuum provided the initial data are of small energy.
Recently, for the case of density-dependent viscosities, \cite{MR4941832} obtained global existence and decay of strong solution to the 3D
compressible NSP equations with arbitrarily large initial data.

On the other hand, the related results with small initial perturbations were well-established.
The global existence of smooth solutions to the compressible NSP system away from vacuum with the optimal time-decay estimates was established in \cite{MR2609958} for the isentropic case,
and in \cite{MR2778703,MR2737817} for non-isentropic case. Then, the pointwise decay estimates of global solutions for bipolar compressible NSP system were obtained by Wu-Wang \cite{MR3687877,MR4113192} for isentropic flow
and non-isentropic flow, respectively. One can find more results concerned with the stability of the wave patterns in \cite{MR3306347,MR3436996,MR4097262,MR4916848,MR4999636} and the references therein.

To the best of our knowledge, however, few results are available on the large-time behavior of global solutions to the 1D non-isentropic compressible Navier-Stokes-Poisson system \eqref{nsp} without imposing smallness conditions or constraint conditions. Even in the case of constant transport coefficients, only Nishida-Smoller type solutions
were known for $\chi=0$, while for $\chi=1$ no results have been available. We recall again that the work \cite{MR3032988} considering system \eqref{nsp} with $\chi=0$ assumed the transport coefficients depend on the density, and also required $\gamma-1$ to be sufficiently small.
See also the work \cite{MR4693518,Zhao-Zhao-2026} considering system \eqref{nsp} with $\chi=1$ asked for the transport coefficients depending on the density and temperature.
Thus, a comprehensive theory for large-data solutions, with uniform bounds and precise large time behaviors, has been lacking.
In this paper, we will give a positive answer to this question for the system \eqref{nsp} with both constant transport coefficients and variable transport coefficients.

\medskip
Now we list the main ideas:
\begin{enumerate}
  \item For the Theorem \ref{main-theorem}, the key step to proving the large time behavior of global solutions is to derive the upper and lower bounds
  of the specific volume $v(t,x)$ which are uniform in time and space.
  A milestone method is that Jiang \cite{MR1671920,MR1912419} constructed a local representation of the specific volume $v(t,x)$ based on
  the global existence result in \cite{MR0468593}. However,
  as \cite{MR3032988} said, ``when the viscosity coefficient $\mu(v)$ is a positive constant,
  the argument in \cite{MR0468593} does not give bounds on $v(t,x)$ because of the nonlocal term $\frac{\Phi_x}{v}$".
  Indeed, we will confront a term $A(t,x)$ as follows:
\begin{equation*}
      A(t,x)=\exp\left\{\frac{1}{\mu}\int_0^t
    \int_{[\widetilde{x}]}^{[\widetilde{x}]+1}\left(\frac{\mu u_x}{v}-p\right)d\xi ds+\frac{1}{\mu}\int_{0}^{t}\int_{x}^{\infty}\frac{\Phi_x}{v}\widetilde{\varphi} d\xi ds\right\},
\end{equation*}
where $\widetilde{\varphi}(x)$ is a cut-off function.
Controlling the additional integral $\int_{0}^{t}\int_{x}^{\infty}\frac{\Phi_x}{v}\widetilde{\varphi} d\xi ds$ from the electric potential is highly nontrivial, as it could potentially grow unbounded as $t\rightarrow\infty$. A possible way to handle this term is an application of H\"{o}lder's inequality, which would require a uniform $L^2$ estimate for $\frac{\Phi_x}{v}$, but the standard energy method only provides control over $\int_{\bbR}\frac{\Phi^2_x}{v}dx$, which is insufficient. Fortunately,
we find a hidden (weighted) $L^2$ estimate for $\frac{\Phi_x}{v}$. More precisely, by carefully examining the Poisson equation $\left(\frac{\Phi_x}{v}\right)_x=1-ve^{-\chi\Phi}$, we discover the following identities:
\begin{align*}
  \int_{\bbR}\frac{\Phi_x^2}{v}dx=
  \begin{cases}
  \int_{\bbR}\left(\frac{\Phi_x}{v}\right)^2dx,\quad &{\rm for}\ \chi=0,\\
  \int_{\bbR}\left(\frac{\Phi_x}{v}\right)^2e^{\Phi}dx+\frac13\int_{\bbR}\frac{\Phi_x^4}{v^3}e^{\Phi}dx,\quad &{\rm for}\ \chi=1.
  \end{cases}
\end{align*}
It is worth noting that we could obtain a uniform upper bound for $|\Phi(t,x)|$ by applying only Corollary \ref{cor_0.5} when $\chi=1$. Hence,
the weighted estimate $\int_{\bbR}\left(\frac{\Phi_x}{v}\right)^2e^{\Phi}dx$ functions the same as $\int_{\bbR}\left(\frac{\Phi_x}{v}\right)^2dx$.
This provides the first uniform control over the potential force.
However, the uniform bound on the potential force is not enough to obtain an exponential decaying estimate for $A(t,x)$.
Actually, since we just regard the potential force as a disturbance term, we have to compare it with the ``effective viscous flux'' $\frac{\mu u_x}{v}-p$
in $A(t,x)$, however, we can not confirm which quantity is larger. This leads to a possible exponential growth on $A(t,x)$. Hence, we require some new ideas.

To obtain the upper bound for $v(t,x)$, we localize the estimate near each point $x$ by combining a spatial antiderivative with averaging.
In Jiang's construction in \cite{MR1671920}, the cutoff function has a transition interval of length one. Here we use the interval $[x-l,x+l]$
and choose a suitable $l>0$, depending only on the initial data. This allows us to use the Poisson equation
\begin{equation*}
  \left(\frac{\Phi_x}{v}\right)_x=1-ve^{-\chi\Phi}.
\end{equation*}
The ingredient is to find a suitable function instead of the cut-off function $\widetilde{\varphi}(x)$ in $A(t,x)$ such that its composition with $\frac{\Phi_x}{v}$ enables integration by parts for the electric potential. With this consideration,
we consider the following operator defined in \eqref{T-2}:
\begin{align*}
   \mathcal{T}_{l,x}h =\frac{1}{2l}\int_{-l}^{l}\left(\int_0^rh(x+z)dz\right)dr
   :=\int_{-l}^lK_l(z)h_z(x+z)dz,
\end{align*}
with
\begin{align*}
  K_l(z)=\begin{cases}
  \frac{(l+z)^2}{4l},&\quad -l\leq z\leq0,\\
   \frac{(l-z)^2}{4l},&\quad 0\leq z\leq l.
  \end{cases}
\end{align*}
Applying this operator to equation $\eqref{nsp}_2$, we obtain a new representation \eqref{v_form2} on the specific volume $v(t,x)$.
What is exciting is that, the function
\begin{align*}
  A_l(t,x)&:=\exp\left\{\frac{1}{\mu}\int_0^tH_l(s,x)ds\right\},\\
   H_l(t,x)&:=\frac{l^2}{6}-R\left\langle\frac{\theta}{v}\right\rangle_{l,x}-\int_{-l}^lK_l(z)(ve^{-\chi\Phi})(t,x+z)dz,
\end{align*}
can be controlled by choosing a suitable $l$ uniformly in $x$. Hence, we obtain the upper bound for the specific volume $v(t,x)$. It is worth noting that this method to derive the upper bound also applies to the Navier-Stokes equations.

Considering the lower bound for $v(t,x)$, we follow Jiang's construction of the local representation in \cite{MR1671920}.
In that work, the uniform lower bound is obtained for an existing global solution by combining a large-time estimate with an exponentially decaying lower bound for the
specific volume. However, due to the lack of suitable space-time estimates for the electric potential force, the above argument cannot be
applied directly to the NSP system. Based on Jiang's construction, it suffices to consider time intervals of a fixed length $h>0$.
More precisely, for $0\leq t\leq h$, the lower bound follows from the initial lower bound of $v_0$; for $t\geq h$, it follows from the
temperature integral over $[t-h,t]$. Thus the same positive lower bound holds on any interval where the local solution exists, with a constant
depending only on the initial data and $h$. Combining this estimate with the local existence result and the a priori estimates,
we can extend the solution successively to obtain a global solution. The lower bound remains unchanged during this
continuity argument and does not depend on the total existence time.


Once the uniform bounds of $v(t,x)$ are obtained, it remains to deduce the uniform bounds from below to above for the temperature $\theta(t,x)$. Although, the proof
can follow the same procedure in \cite{MR3466845}, the electric potential $\Phi(t,x)$ brings major problems.
For example, the following terms
\begin{equation*}
  \int_0^T\int_{\bbR}\frac{\Phi_x}{v}u(\theta-2)_+,\quad \int_0^t\int_{\bbR}\frac{\Phi_x}{v}u^3
\end{equation*}
in \eqref{J-6} and \eqref{pro2.4.9}, respectively, should be treated carefully.

When $\chi=0$, only by using equation $\eqref{nsp}_4$ and integration by parts,
one can not handle this term, due to the lack of $L^2$ estimates in temporal and spatial for $\Phi_x(t,x)$.
To resolve these difficulties, our methods are highly dependent on the Lemma \ref{lem-2.2}, the domain decomposition technique.
When $\chi=1$, the Boltzmann relation has brought both benefits and drawbacks. More precisely, the Boltzmann relation causes Lemma \ref{lem-2.2} to fail.
Instead, a damping effect is induced which provides a good estimate for $\Phi_x(t,x)$.
Moreover, a new dissipative term
\begin{align*}
  \int_0^T\int_{\bbR}(1-ve^{-\chi\Phi})\ln(ve^{-\chi\Phi})\leq 0
\end{align*}
is derived, which provide the large time behavior of $\|1-ve^{-\chi\Phi}\|_{L^2}$.
It is worth noting that similar dissipative term has been derived when we study the isentropic case in \cite{MR4458894,MR4755045}. The detail is proposed in Section 2.
  \item When the transport coefficients $\mu(v,\theta)$ and $\kappa(v,\theta)$ depend on the specific volume $v(t,x)$ and temperature $\theta(t,x)$, the local representation for the specific volume $v(t,x)$
is not valid any more. Instead, the strategy introduced by Kanel$'$ \cite{MR0227619} which has been used frequently (c.f.\cite{MR3032988,MR3861291,MR3225502}) to derive the upper and lower bounds for the
specific volume is applied. In contrast to the case of constant transport coefficients, we obtain more involved estimates for $\Phi(t,x)$ in \eqref{vx-v} as follows:
\begin{equation*}
  \int_0^T\int_{\mathbb{R}}\chi\frac{\mu(e^{\chi\Phi},\theta)}{v}\Phi^2_x,\quad
  \int_0^T\int_{\mathbb{R}}\left(g(v,\theta)-g(e^{\chi\Phi},\theta)\right)\left(ve^{-\chi\Phi}-1\right)\geq0,
\end{equation*}
which are crucial to close the energy estimates.

The most difficult problem lies in proving the upper bound for the temperature due to the nonlocal terms such as \eqref{3.4-10} and \eqref{3.4-21}. Compared with
Theorem \ref{main-theorem} whose upper and lower bounds of the specific volume are independent of the space and time, the temperature-dependent bounds of the specific volume in Lemma \ref{lem-3.3}
ask for more precise estimates. Thanks to the Lemma \ref{lem-2.2}, we can
transform the nonlocal terms into several terms that are relatively easy to handle.
Since we have not yet obtained the uniform bounds for the specific volume,
the term $\mathcal{K}_{3}$ in \eqref{3.4-13} as well as the term $\mathcal{L}_{1}$ in \eqref{3.4-L1} should be handled more carefully compared with
the estimates \eqref{pro2.4.5} and \eqref{K-1}.
We find that by dividing the integrant into two parts for $v\geq 2$ and $v\leq 2$, these difficult terms can be successfully controlled by an explicit expression for $\|\theta\|_{L^{\infty}_{t,x}}$ denoted by $F(\theta)$ in \eqref{F-theta}
and finally derive the $L^2-$type energy estimates for the
temperature $\theta(t,x)$.  The main procedure is carried out in Section 3.
\end{enumerate}

\subsection{Outline of the paper and notations}
The rest of the paper is organized
as follows: Section 2 focuses on the Cauchy problem \eqref{nsp} and \eqref{initial-data} with constant transport coefficients and is devoted to the proof of Theorem \ref{main-theorem}.
A physical case for transport coefficients depending on specific volume and temperature is considered in section 3.
In the appendix, we prepare some Lemmas which are crucial to prove the uniform bounds for $v(t,x)$, $\theta(t,x)$ and $\Phi(t,x)$.
\\

 \textit{Notations.} For notational simplicity, throughout this paper, $N$, $C$, $m_i$, $C_{i}$ and $B_i(i\in\bbR)$ denote some generic positive constants changing from line to line.
$C(*)$ represents that the positive constant $C$ depends on $*$. $A\lesssim B$ means that there is a generic
constant $C>0$ such that $A\leq CB$. $A\thicksim B$ means $A\lesssim B$ and $B\lesssim A$.
$L^q(\mathbb{R})(1\leq q\leq \infty)$ stands for the usual Lebesgue space with norm $\|\cdot\|_{L^q}$,
and $H^k(\mathbb{R})(k\in \mathbb{N})$ is the usual Sobolev space in the $L^2$ sense with norm $\|\cdot\|_{k}$. $L^{\infty}(\bbR)$ and $L^{\infty}_{t,x}(\bbR)$ means the Lebesgue space with respect to $x$ and both of $t$ and $x$, respectively.

\section{The proof of Theorem \ref{main-theorem}}
In this section, we will prove the Theorem \ref{main-theorem}. The local existence results and the a priori estimates will be presented. It is crucial to derive the uniform upper and lower bounds
for the specific volume $v(t,x)$ and temperature $\theta(t,x)$ independently of space and time which contributes to deriving the large time behavior of the global solutions.

\begin{lemma}[c.f. \cite{MR3032988}]\label{lem-2.1}
Under the assumptions in Theorem \ref{main-theorem}, there exists a small positive contant $T$ depends only on the initial data,
such that the Cauchy problem \eqref{nsp}-\eqref{initial-data} admits a unique smooth solution $(v(t,x),u(t,x),\theta(t,x),\Phi(t,x))$ which satisfies
\begin{itemize}
  \item When $\chi=0$,
   \begin{equation}
  \begin{cases}
  (v-1,u,\theta-1)\in C^0(0,T;H^1(\bbR)),\quad \Phi_{x}\in C^0(0,T;H^1(\bbR)),\\
  v_x\in L^2(0,T;L^2(\bbR)),\quad (u_x,\theta_x)\in L^2(0,T;H^1(\bbR)).
  \end{cases}
\end{equation}
  \item When $\chi=1$,
   \begin{equation}
  \begin{cases}
  (v-1,u,\theta-1)\in C^0(0,T;H^1(\bbR)),\quad \Phi\in C^0(0,T;H^2(\bbR)),\\
  v_x\in L^2(0,T;L^2(\bbR)),\quad (u_x,\theta_x,\Phi_x)\in L^2(0,T;H^1(\bbR)).
  \end{cases}
\end{equation}
\end{itemize}
\end{lemma}


\subsection{The a priori estiamtes}
Now, we are concerned with the basic energy estimates. We begin with a convex entropy defined as follows:
\begin{equation}\label{entropy}
  \eta(v,u,\theta)=R\phi(v)+\frac{u^2}{2}+C_v\phi(\theta)\quad  \mathrm{with} \quad \phi(x)=x-\ln x-1.
\end{equation}
Direct calculation shows that
\begin{equation}\label{eta}
  \eta_t+\left\{\left(\frac{R\theta}{v}-R\right)u\right\}_x-\left(\frac{\mu uu_x}{v}+\frac{\kappa(\theta-1)\theta_x}{v\theta}\right)_x+\left(\frac{\mu u_x^2}{v\theta}+\frac{\kappa\theta_x^2}{v\theta^2}\right)=\frac{u\Phi_x}{v}.
\end{equation}
Besides, by integration by parts, we have
\begin{align}
  \int_{\bbR}\frac{u\Phi_x}{v}dx=&\int_{\bbR} u\frac{\Phi_x}{v}\left(\left(\frac{\Phi_x}{v}\right)_x+ve^{-\chi\Phi}\right)dx
  =\int_{\bbR} u\Phi_x e^{-\chi\Phi} dx-\frac12\int u_x\left(\frac{\Phi_x}{v}\right)^2dx\nonumber\\
  =&\frac{d}{dt}\int_{\mathbb{R}}\left\{-\chi v\left(1- e^{-\Phi}-\Phi e^{-\Phi}\right)-\frac{\Phi_x^2}{2v}\right\} dx.\label{u-phi-first}
\end{align}
The detailed derivation of \eqref{u-phi-first} can be found in our works \cite{MR4693518,MR4755045}.



Combining \eqref{eta} and \eqref{u-phi-first}, the following corollary can be derived by integrating with respect to $t$ and $x$ over $[0,T]\times \bbR$.
\begin{corollary}\label{cor_0.5}
For each time $t\in[0,T]$, the following estimate holds:
  \begin{align}
  &\sup_{0\leq t\leq T}\int_{\bbR}\left[\eta(v,u,\theta)+\frac{\Phi_x^2}{2v}+\chi v(1-\Phi e^{-\Phi}-e^{-\Phi})\right](t,x)dx
    +\int_0^{t}\int_{\bbR}\left[\frac{\mu u_{x}^2}{v\theta}
  +\frac{\kappa\theta_x^2}{v\theta^2}\right]dxds \nonumber\\
  \leq &\int_{\bbR}\left[\eta(v,u,\theta)+\frac{\Phi_x^2}{2v}+\chi v(1-\Phi e^{-\Phi}-e^{-\Phi})\right](0,x)dx:=e_1.\label{cor-lag-1}
\end{align}
\end{corollary}

However, the estimates \eqref{cor-lag-1} are not good enough to derive the uniform upper and lower bounds for the specific volume $v(t,x)$
in Lemma \ref{lem_boun}. Fortunately, we find a hidden relation below by applying equation $\eqref{nsp}_4$ which gives a favorable $L^2$ estimates for $\frac{\Phi_x}{v}$:
\begin{itemize}
  \item When $\chi=0$,
  \begin{align}\label{new-relation-0}
  \int_{\bbR} \frac{\Phi_x^2}{v}dx-\int_{\bbR}\left(\frac{\Phi_x}{v}\right)^2dx&=\int_{\bbR}\left(\frac{\Phi_x}{v}\right)^2(v-1)dx\notag\\
  &=-\int_{\bbR}\left(\frac{\Phi_x}{v}\right)^2\left(\frac{\Phi_x}{v}\right)_xdx=0.
\end{align}
  \item When $\chi=1$,
  \begin{align}\label{new-relation}
  \int_{\bbR}\frac{\Phi_x^2}{v}dx-\int_{\bbR}\left(\frac{\Phi_x}{v}\right)^2e^{\Phi}dx&=\int_{\bbR}\left(\frac{\Phi_x}{v}\right)^2e^{\Phi}(ve^{-\Phi}-1)dx\notag\\
  &=-\int_{\bbR}\left(\frac{\Phi_x}{v}\right)^2e^{\Phi}\left(\frac{\Phi_x}{v}\right)_x dx=\frac13\int_{\bbR}\frac{\Phi_x^4}{v^3}e^{\Phi}dx\geq0.
\end{align}
\end{itemize}
Notice that the weighted estimates of potential in \eqref{new-relation} are positive.
Since a uniform upper bound for $|\Phi(t,x)|$ can be derived by only applying Corollary \ref{cor_0.5},
the term $\int_{\bbR}\left(\frac{\Phi_x}{v}\right)^2e^{\Phi}dx$ functions the same as $\int_{\bbR}\left(\frac{\Phi_x}{v}\right)^2dx$.\\

Using the identities \eqref{new-relation-0} and \eqref{new-relation} as well as Lemma \ref{lem2.2}, we can obtain a new corollary below:
\begin{corollary}\label{cor_0}
For each time $t\in[0,T]$, the following estimate holds:
\begin{itemize}
  \item When $\chi=0$, we have
  \begin{align}
  &\sup_{0\leq t\leq T}\int_{\bbR}\left[\eta(v,u,\theta)+\frac12\left(\frac{\Phi_x}{v}\right)^2\right](t,x)dx
    +\int_0^{t}\int_{\bbR}\left[\frac{\mu u_{x}^2}{v \theta}
  +\frac{\kappa\theta_x^2}{v\theta^2}\right]dxds \nonumber\\
  \leq &\int_{\bbR}\left[\eta(v,u,\theta)+\frac12\left(\frac{\Phi_x}{v}\right)^2\right](0,x)dx:=e_2.\label{vol_0}
\end{align}
  \item When $\chi=1$, using Lemma \ref{lem2.2}, we obtain
\begin{align}
  &\sup_{0\leq t\leq T}\int_{\bbR}\left[\eta(v,u,\theta)+\frac12\left(\frac{\Phi_x}{v}\right)^2e^{\Phi}+\frac16\frac{\Phi_x^4}{v^3}e^{\Phi}+v(1-\Phi e^{-\Phi}-e^{-\Phi})\right](t,x)dx \nonumber\\
    &\hspace{1cm}+\int_0^{t}\int_{\bbR}\left[\frac{\mu u_{x}^2}{v\theta}
  +\frac{\kappa\theta_x^2}{v\theta^2}\right]dxds \nonumber\\
  \leq &C(\underline{V},\overline{V})\left[\int_{\bbR}\eta(v,u,\theta)(0,x)dx+\|v_0-1\|_{L^2}^2+\|v_0-1\|_{L^2}^6\right]:=e_3.\label{cor-lag-2}
\end{align}
\end{itemize}
\end{corollary}

For any $\widetilde{x}\in\bbR$, it can be derived from Corollary $\ref{cor_0}$ that
\begin{equation}\label{vol_1}
  \int_{[\widetilde{x}]}^{[\widetilde{x}]+1}\phi(v)(t,x)dx,\quad \int_{[\widetilde{x}]}^{[\widetilde{x}]+1}\phi(\theta)(t,x)dx\leq \frac{\max\{e_2,e_3\}}{\min\{R,c_v\}}:=M_*.
\end{equation}
Moreover, utilizing $\eqref{vol_1}$ and applying Jensen's inequality
to the convex function $\phi$, we obtain:
\begin{equation*}
  \phi\left(\int_{[\widetilde{x}]}^{[\widetilde{x}]+1}v(t,x)dx\right),\quad \phi\left(\int_{[\widetilde{x}]}^{[\widetilde{x}]+1}\theta(t,x)dx\right)\leq M_*,
\end{equation*}
which gives $\alpha_1\leq \int_{[\widetilde{x}]}^{[\widetilde{x}]+1}v(x,t)dx,\quad \int_{[\widetilde{x}]}^{[\widetilde{x}]+1}\theta(x,t)dx \leq \alpha_2$, where
$\alpha_1,\alpha_2$ are two positive roots of the equation $\phi(x)=M_*$. Since $\phi(1)=0$, we can obtain $0<\alpha_1<1<\alpha_2$. Hence by terms of the mean value theorem, for
each $t\geq0$ there are points $a_{[\widetilde{x}]}(t),b_{[\widetilde{x}]}(t)\in[[\widetilde{x}],[\widetilde{x}]+1]$ such that
\begin{equation*}
  \alpha_1\leq v(t,a_{[\widetilde{x}]}(t)),\quad \theta(t,b_{[\widetilde{x}]}(t))\leq \alpha_2.
\end{equation*}

From above analysis, we thus have proved
\begin{lemma}\label{2.3.0}
Let $\alpha_1,\alpha_2$ be two positive roots of the equation $\phi(x)=M_*$. Then, for each $t\in[0,T]$,
\begin{equation}\label{vol_2}
  \alpha_1\leq \int_{[\widetilde{x}]}^{[\widetilde{x}]+1}v(t,x)dx,\quad \int_{[\widetilde{x}]}^{[\widetilde{x}]+1}\theta(t,x)dx\leq \alpha_2,\quad t\geq0,
\end{equation}
and there are points $a_{[\widetilde{x}]}(t),b_{[\widetilde{x}]}(t)\in[[\widetilde{x}],[\widetilde{x}]+1]$ such that
\begin{equation}\label{vol_3}
  \alpha_1\leq v(t,a_{[\widetilde{x}]}(t)),\quad \theta(t,b_{[\widetilde{x}]}(t))\leq \alpha_2,\quad t\geq0.
\end{equation}
\end{lemma}

Next, we deduce two different local representations of $v(t,x)$ which are crucial to proving the
upper and lower bounds for this quantity, respectively:

{\it 1. For the first one}, we will apply Jiang's argument in  \cite{MR1671920,MR1912419} for
fixed domain. This local representation of $v(t,x)$ is used to deriving the lower bound.
To this end, we introduce the cutoff function
  $\varphi\in W^{1,\infty}(\mathbb{R})$ with
parameter $x\in\mathbb{R}$ by
\begin{equation}\label{varphi}
  \varphi(x)=
  \begin{cases}
    1,\quad& x\leq[\widetilde{x}]+1,\\
    [\widetilde{x}]+2-x,\quad& [\widetilde{x}]+1\leq x\leq[\widetilde{x}]+2,\\
    0,\quad& x\geq [\widetilde{x}]+2.
  \end{cases}
\end{equation}
where $[x]$ denotes the largest integer that is less or equal to $x$. For simplicity we denote $\Omega_{[\widetilde{x}]}:=([\widetilde{x}]-1,[\widetilde{x}]+1)$, then we have the following lemma:
\begin{lemma}
We have
\begin{equation}\label{v_form}
  v(t,x)=Y_s(t)M_s(t,x)v(s,x)+\frac{R}{\mu}\int_s^t
  Y_r(t)M_r(t,x)\theta(r,x)dr
\end{equation}
for all $0\leq s \leq t\leq T$ and $x\in\overline{\Omega}_{[\widetilde{x}]}$, where
  \begin{align*}
    Y_s(t)&:=\exp\left\{\frac{1}{\mu}\int_s^t\int_{[\widetilde{x}]+1}^{[\widetilde{x}]+2}\left(\mu\frac{u_y}{v}-p\right)(\tau,y)dyd\tau\right\},\\
    M_s(t,x)&:=\exp\left\{-\frac{1}{\mu}\left(\int_x^{\infty}u(t,y)\varphi(y)dy-\int_x^{\infty}u(s,y)\varphi(y)dy\right)+\frac{1}{\mu}\int_s^t\int_x^{\infty}\frac{\Phi_x}{v}(t,y)\varphi(y)dyd\tau\right\}.
  \end{align*}
\end{lemma}

\begin{proof}
We  multiply $\eqref{nsp}_2$ by $\varphi$ to get
 \begin{equation} \label{ident1}
   (\varphi u)_t=\left[\left(\mu\frac{u_x}{v}-p\right)\varphi\right]_x
   -\varphi_x\left(\mu\frac{u_x}{v}-p\right)+\frac{\Phi_x}{v}\varphi.
 \end{equation}
In view of the definition of $\varphi(x)$ and \eqref{nsp}$_1$,
we integrate \eqref{ident1} over
$[x,\infty)$ to get
 \begin{align*}
   (\ln v(t,x))_t&=\frac{R}{\mu}\frac{\theta}{v}+\frac{1}{\mu}\left(\int_{[\widetilde{x}]+1}^{[\widetilde{x}]+2}\left(\mu\frac{u_y}{v}-p\right)(t,y)dy+\int_x^{\infty}\frac{\Phi_x}{v}(t,y)\varphi(y)dy-\left(\int_x^{\infty}u(t,y)\varphi(y)dy\right)_t\right)\\
   &=\frac{R}{\mu}\frac{\theta}{v}+(\ln (Y_s(t)M_s(t,x)))_t.
 \end{align*}
Rewriting the above identity, we have
\begin{align}\label{v-1}
  \left(\frac{v(t,x)}{Y_s(t)M_s(t,x)}\right)_t=\frac{R}{\mu}\frac{\theta(t,x)}{Y_s(t)M_s(t,x)}.
\end{align}
Integrating \eqref{v-1} with respect to $t$ over $[s,t]$ yields
the desired representation \eqref{v_form}, this completes the proof of the lemma.
\end{proof}

{\it 2. For the second one}, it is applied to obtain the upper bound of $v(t,x)$ by
localizing near each point $x$ with a suitable scale $l>0$. Note that the identity \eqref{v_form}
is invalid to derive the upper bound, some new ideas are needed. Before giving the representation, we list the following
notation, for any smooth function $h(t,x)$ and fixed $l>0$ which will be chosen later, we define
\begin{align}\label{T}
   \langle h \rangle_{l,x}(x):=\frac{1}{2l}\int_{x-l}^{x+l}h(y)dy,\quad \mathcal{T}_{l,x}h(x):=\frac{1}{2l}\int_{x-l}^{x+l}\left(\int_x^yh(z)dz\right)dy.
\end{align}
Fixing $x$ and letting $r=y-x$, we obtain
\begin{align}\label{T-2}
   \mathcal{T}_{l,x}h =&\,\frac{1}{2l}\int_{-l}^{l}\left(\int_0^rh(x+z)dz\right)dr\notag\\
   =&\,\frac{1}{2l}\int_{0}^{l}\left(\int_0^rh(x+z)dz\right)dr-\frac{1}{2l}\int_{-l}^{0}\left(\int_r^0h(x+z)dz\right)dr\notag\\
   =&\,\frac{1}{2l}\int_{0}^{l}(l-z)h(x+z)dz-\frac{1}{2l}\int_{-l}^0(l+z)h(x+z)dz\notag\\
   =&\,-\int_{0}^{l}\left(\frac{(l-z)^2}{4l}\right)_zh(x+z)dz-\int_{-l}^0\left(\frac{(l+z)^2}{4l}\right)_zh(x+z)dz\notag\\
   =&\,-\left.\left[\frac{(l-z)^2}{4l}h(x+z)\right]\right|_{0}^l-\left.\left[\frac{(l+z)^2}{4l}h(x+z)\right]\right|_{-l}^0\notag\\
   &\,+\int_{0}^{l}\frac{(l-z)^2}{4l}h_z(x+z)dz+\int_{-l}^0\frac{(l+z)^2}{4l}h_z(x+z)dz\notag\\
   :=&\,\int_{-l}^lK_l(z)h_z(x+z)dz,
\end{align}
with
\begin{align*}
  K_l(z)=\begin{cases}
  \frac{(l+z)^2}{4l},&\quad -l\leq z\leq0,\\
   \frac{(l-z)^2}{4l},&\quad 0\leq z\leq l.
  \end{cases}
\end{align*}
We have used the fact that $K_l(-l)=K_l(l)=0$ and $h(x)$ is a continue function. Notice that the kernel $K_l(z)$ satisfies
\begin{align}\label{K-property}
   \int_{-l}^lK_l(z)dz=\frac{l^2}{6},\quad \int_{-l}^0(K'_l)^2(z)dz+\int_{0}^l(K'_l)^2(z)dz=\frac{l}{6}.
\end{align}

With the above preparation, we have the following lemma:
\begin{lemma}\label{lem-v_form2}
We have
\begin{equation}\label{v_form2}
  v(t,x)=B_l(t,x)A_l(t,x)+\frac{R}{\mu}\int_0^t\frac{B_l(t,x)A_l(t,x)}{B_l(s,x)A_l(s,x)}\theta(s,x)ds
\end{equation}
for all $t\in[0,T]$ and $x\in \bbR$, where
\begin{align}
   B_l(t,x)&:=v_0(x)\exp\left\{\langle\ln v-\ln v_0\rangle_{l,x}-\frac{1}{\mu}\mathcal{T}_{l,x}(u-u_0)\right\},\\
   A_l(t,x)&:=\exp\left\{\frac{1}{\mu}\int_0^tH_l(s,x)ds\right\},\\
   H_l(t,x)&:=\frac{l^2}{6}-R\left\langle\frac{\theta}{v}\right\rangle_{l,x}-\int_{-l}^lK_l(z)(ve^{-\chi\Phi})(t,x+z)dz
\end{align}
\end{lemma}
\begin{proof}
Applying $\mathcal{T}_{l,x}$ to $\eqref{nsp}_2$, we obtain
\begin{align*}
   (\mathcal{T}_{l,x}u)_t=\langle\sigma\rangle_{l,x}-\sigma(t,x)+\mathcal{T}_{l,x}\frac{\Phi_x}{v},
\end{align*}
with $\sigma=\frac{\mu u_x}{v}-p$. Note that $\langle\sigma\rangle_{l,x}=\mu (\langle\ln v\rangle_{l,x})_t-R\left\langle\frac{\theta}{v}\right\rangle_{l,x}$, we have
\begin{align}\label{v-form-1}
   (\ln v)_t&=\frac{R}{\mu}\frac{\theta}{v}+(\langle\ln v\rangle_{l,x})_t-\frac{1}{\mu}(\mathcal{T}_{l,x}u)_t+\frac{1}{\mu}H_l(t,x)\notag\\
   &=\frac{R}{\mu}\frac{\theta}{v}+(\ln(B_lA_l))_t,
\end{align}
where
\begin{align*}
   H_l(t,x)&=-R\left\langle\frac{\theta}{v}\right\rangle_{l,x}+\int_{-l}^{l}K_l(z)\left(\frac{\Phi_z}{v}\right)_z(t,x+z)dz,\\
   &=\frac{l^2}{6}-R\left\langle\frac{\theta}{v}\right\rangle_{l,x}-\int_{-l}^lK_l(z)(ve^{-\chi\Phi})(t,x+z)dz.
\end{align*}
Rewriting \eqref{v-form-1}, we have
\begin{equation*}
  \left(\frac{v(t,x)}{B_lA_l}\right)_t=\frac{R}{\mu}\frac{\theta(t,x)}{B_lA_l}.
\end{equation*}
Integrating the above identity with respect to $t$ over $[0,t]$ yields
the desired representation \eqref{v_form2}, this completes the proof of the lemma.
\end{proof}

Now in order to bound $v(t,x)$ pointwisely, we first show the exponential decay of $A(t)$,
then we use the representation \eqref{v_form2} to obtain the uniform upper bound on $v(t,x)$. The uniform lower bound on $v(t,x)$ depends on
the representation \eqref{v_form}.
\begin{lemma} \label{lem_boun}
For any fixed $h>0$, there exist two positive constants $\underline{v}_h$ and $\overline{v}$ independent of time $T$ such that
     \begin{equation}\label{boun_v}
   \underline{v}_h\leq v(t,x)\leq \overline{v}\quad {\rm for\ all}\ (t,x)\in[0,T]\times\bbR.
  \end{equation}
\end{lemma}
\begin{proof}
The proof is divided into three steps.

\noindent{\bf Step 1}. Firstly, we derive the uniform upper bound for $v(t,x)$.
The operator $\mathcal{T}_{l,x}$ is a bounded operator such that
\begin{equation*}
  |\mathcal{T}_{l,x}u|\leq \sqrt{\frac{l}{6}}\|u\|_{L^2}\leq C\sqrt{\frac{l}{6}}.
\end{equation*}
Here, $C$ is a generic constant depends only on the initial perturbation. Then, it is easy to calculate that
\begin{align*}
  |\ln z|\leq \phi(z)+z+1\leq C(\phi(z)+1),
\end{align*}
so we have
\begin{align*}
  |\langle\ln v\rangle_{l,x}|\leq \frac{1}{2l}\int_{x-l}^{x+l}|\ln v(t,y)| dy\leq C+\frac{C}{l}\int_{\bbR}\phi(v)dy\leq C(l).
\end{align*}
Hence, for the term $B_l(t,x)$, we have
\begin{equation}\label{est_B}
  \underline{V}e^{-C(l)}\leq B_l(t,x)\leq \overline{V}e^{C(l)} \quad {\rm for\ all}\ (t,x)\in[0,T]\times\bbR.
\end{equation}
Here, $C(l)$ depends only on $l$ and the initial perturbation.

Next, we will derive a lower bound for $\theta(t,x)$. To achieve this goal,
we divide the time $t\in[0,T]$ into two parts according to the dissipation in the following way:
\begin{align}\label{good time}
  \mathcal{G}(t):=\left\{t\in[0,T]: \int_{\bbR}\frac{\kappa\theta_x^2}{v\theta^2}dx \leq 1\right\},\quad \mathcal{B}(t):=\left\{t\in[0,T]: \int_{\bbR}\frac{\kappa\theta_x^2}{v\theta^2}dx > 1\right\}.
\end{align}
For any $x\in[[\widetilde{x}],[\widetilde{x}]+1]$, we have
\begin{align}\label{theta-lower}
  |\ln \theta(t,x)-\ln\theta(t,b_{[\widetilde{x}]}(t))|\leq \left(\int_{[\widetilde{x}]}^{[\widetilde{x}]+1}\frac{\theta_x^2}{v\theta^2}dy\right)^{\frac12}\left(\int_{[\widetilde{x}]}^{[\widetilde{x}]+1}vdy\right)^{\frac12}
  \leq \sqrt{\frac{\alpha_2}{\kappa}}\left(\int_{\bbR}\frac{\kappa\theta_x^2}{v\theta^2}dx\right)^{\frac12}.
\end{align}
In the good time region $\mathcal{G}(t)$, it follows \eqref{theta-lower} that
\begin{equation}\label{low-theta-good}
  \ln \theta(t,x)\geq \ln \alpha_1-\sqrt{\frac{\alpha_2}{\kappa}}\Rightarrow \theta(t,x)\geq \alpha_1e^{-\sqrt{\frac{\alpha_2}{\kappa}}},\quad \forall x\in\bbR.
\end{equation}
On the other hand, the measure of the bad set $\mathcal{B}(t)\cap [0,T]$ admits a uniform upper bound, indeed,
\begin{align}\label{bad-est}
  |\mathcal{B}\cap [0,T]|=\int_0^T 1_{\mathcal{B}}ds< \int_0^T\int_{\bbR}\frac{\kappa\theta_x^2}{v\theta^2}dx\leq \max\{e_2,e_3\}.
\end{align}

With the lower bound \eqref{low-theta-good} on $\theta(t,x)$, in the good time region $\mathcal{G}(t)$, we have
\begin{align*}
  &\,R\left\langle\frac{\theta}{v}\right\rangle_{l,x}+\int_{-l}^lK_l(z)(ve^{-\chi\Phi})(t,x+z)dz\\
  \geq&\, \int_{-l}^{l}\left[\frac{R\alpha_1e^{-\sqrt{\frac{\alpha_2}{\kappa}}}}{2lv(t,x+z)}+m_{\Phi}K_l(z)v(t,x+z)\right]dz\\
  \geq&\, \int_{-l}^{l}2\sqrt{R\alpha_1 e^{-\sqrt{\frac{\alpha_2}{\kappa}}}m_{\Phi}}\sqrt{\frac{K_l(z)}{2l}}dz\\
  \geq&\, l\sqrt{\frac{R\alpha_1 e^{-\sqrt{\frac{\alpha_2}{\kappa}}}m_{\Phi}}{2}},
\end{align*}
where we have used $\int_{-l}^l\sqrt{\frac{K_l(z)}{2l}}dz=\frac{l}{2\sqrt{2}}$, and the inequality $x+y\geq 2\sqrt{xy}$ for $x,y\geq 0$.
Applying Lemma \ref{lem2.3} yields $e^{-\chi\Phi}\geq m_{\Phi}:=\min\{e^{-\overline{\Phi}},1\}$ for different values of $\chi$.

Now, we can choose $l:=\min\{1,\sqrt{R\alpha_1e^{-\sqrt{\frac{\alpha_2}{\kappa}}}m_{\Phi}}\}$, so in the good time region $\mathcal{G}$, the quantity
\begin{align*}
  H_l(t,x)\leq \frac{l^2}{6}-l\sqrt{\frac{R\alpha_1 e^{-\sqrt{\frac{\alpha_2}{\kappa}}}m_{\Phi}}{2}}\leq \frac{l^2}{6}-\frac{l^2}{\sqrt{2}}\leq -\frac{l^2}{2}.
\end{align*}
In the bad region $\mathcal{B}$, apart from $\frac{l^2}{6}$, the other terms are non-positive, we omit them to get
\begin{align*}
  H_l(t,x)\leq \frac{l^2}{6}.
\end{align*}
Combining the above estimates, for any $0\leq s\leq t\leq T$, we obtain from \eqref{bad-est} that
\begin{align}\label{A-l}
  \int_s^tH_l(r,x)dr&=\int_{\mathcal{G}\cap[s,t]}H_l(r,x)dt+\int_{\mathcal{B}\cap[s,t]}H_l(r,x)dr\notag\\
  &\leq -\frac{l^2}{2}(t-s-|\mathcal{B}\cap[s,t]|)+|\mathcal{B}\cap[s,t]|\frac{l^2}{6}\notag\\
  &\leq -\frac{l^2}{2}(t-s)+\frac23l^2\max\{e_2,e_3\},
\end{align}
which yields
\begin{equation}\label{pro7}
    0\leq A({t,x})\leq C\mathrm{e}^{-Ct},\quad
    \frac{A(t,x)}{A(s,x)}\leq C\mathrm{e}^{-C(t-s)}
    \qquad {\rm for\ all}\ 0\leq s\leq {t}\leq T.
  \end{equation}

 \noindent{\bf Step 2}.
    Plugging \eqref{est_B} and \eqref{pro7} into \eqref{v_form2}, we infer
\begin{align}\label{pro2.2.30}
     &C\int_0^t\frac{A(t,x)}{A(s,x)}\theta(s,x)ds\leq v(t,x)\leq C+C\int_0^t\theta(s,x)\mathrm{e}^{-C(t-s)}ds
\end{align}
for all $(t,x)\in [0,T]\times\bbR$.
In light of the fundamental theorem of calculus, we infer
from \eqref{vol_2} that for any
 $x\in[[\widetilde{x}],[\widetilde{x}]+1]$,
\begin{align}
    \left|\theta(s,x)^{\frac{1}{2}}-\theta(s,b_{[\widetilde{x}]}(s))^{\frac{1}{2}}
    \right|
    \lesssim~&
     \int_{[\widetilde{x}]}^{[\widetilde{x}]+1}\theta^{-\frac12}|\theta_{\xi}(s,\xi)|d\xi\nonumber\\
    \lesssim~&
    \left[\int_{[\widetilde{x}]}^{[\widetilde{x}]+1}\frac{\theta_{\xi}^2}{v\theta^2}(s,\xi)d\xi\right]^{\frac{1}{2}}
    \left[\int_{[\widetilde{x}]}^{[\widetilde{x}]+1}{v\theta}(s,\xi)d\xi\right]^{\frac{1}{2}}\nonumber\\
    \lesssim ~& \alpha_2^{\frac12}
    \sup_{[[\widetilde{x}],[\widetilde{x}]+1]}v^{\frac12}(s,\cdot)
    \left[\int_{[\widetilde{x}]}^{[\widetilde{x}]+1}\frac{\theta_{\xi}^2}{v\theta^2}(s,\xi)d\xi\right]^{\frac{1}{2}},\label{pro14}
\end{align}
  where we have used
  $b_{[\widetilde{x}]}(s)\in [[\widetilde{x}],[\widetilde{x}]+1]$.
   Combining \eqref{pro14} with \eqref{vol_3} gives
\begin{equation}\label{pro2.2.305}
    \theta(s,x)\leq~
    2\alpha_2+2\alpha_2\sup_{[[\widetilde{x}],[\widetilde{x}]+1]}v(s,\cdot)
    \int_{[\widetilde{x}]}^{[\widetilde{x}]+1}\frac{\theta_{\xi}^2}{v\theta^2}(s,\xi)d\xi,\quad x\in[[\widetilde{x}],[\widetilde{x}]+1].
  \end{equation}
Hence, gathering all the intervals, we obtain
\begin{equation}\label{pro2.2.31}
    \theta(s,x)\leq~
    2\alpha_2+2\alpha_2\sup_{\bbR}v(s,\cdot)
    \int_{\bbR}\frac{\theta_{\xi}^2}{v\theta^2}(s,\xi)d\xi,\quad x\in\bbR.
  \end{equation}

 We substitute \eqref{pro2.2.31} into \eqref{pro2.2.30} to obtain
 \begin{equation} \label{pro13}
  \sup_{\bbR}v(t,\cdot)\lesssim C+C\int_0^t
  \sup_{\bbR}v(s,\cdot)\int_{\bbR}
    \frac{\theta_{\xi}^2}{v\theta^2}(s,\xi)d\xi ds.
 \end{equation}
 Applying Gronwall's inequality to \eqref{pro13},
 we can deduce from \eqref{vol_0} and \eqref{cor-lag-2}
 that
 \begin{equation}\label{upper1}
  \sup_{\bbR} v (t,\cdot)\leq C, \quad {\rm for\ all}\ t\in[0,T],
 \end{equation}
 where
 $C>0$ is a generic constant independent of $t$, $T$.

\noindent{\bf Step 3}. We will derive the lower bound for the specific volume $v(t,x)$.
Note that the representation \eqref{v_form2} is not useful to derive the lower bound, instead, we will apply the
representation \eqref{v_form} for $0\leq s\leq t\leq T$ with $t-s\leq h$ with fixed $h>0$. Different from \eqref{v_form2}, we separate the electric potential force with the pressure,
and define two new functions $M_s(t,x)$ and $Y_s(t,x)$.

For $M_s(t,x)$, for any $x\in \Omega_{[\widetilde{x}]}$, it derives from the H\"{o}lder's inequality, \eqref{vol_0} and \eqref{cor-lag-2} that
\begin{align*}
  |J(t,x)|\leq C\|u\|_{L^2}\leq C_J,\quad |G(t,x)|\leq C\left\|\frac{\Phi_x}{v}\right\|_{L^2}\leq C_G.
\end{align*}
When $0\leq t-s\leq h$, we get
\begin{align}\label{M-s}
  e^{-C(h)}\leq M_s(t,x)\leq e^{C(h)}.
\end{align}
Here $C(h)=\frac{2C_J+hC_G}{\mu}>0$ is a generic constant which depends only on the initial perturbation, and are independent of $[\widetilde{x}]$, $t$ and $T$.

For $Y_s(t,x)$, we rewrite \eqref{v_form} as follows:
\begin{equation*}
  \frac{1}{Y_s(t)}\frac{v(t,x)}{M_s(t,x)}=v(s,x)+\frac{R}{\mu}\int_s^t
  \frac{1}{Y_s(r)}\frac{\theta(r,x)}{M_s(r,x)}dr
\end{equation*}
Integrating over $\overline{\Omega}_{[\widetilde{x}]}$, it can be derived from \eqref{M-s} and \eqref{vol_2} that,
\begin{align*}
  \frac{1}{Y_s(t)}\alpha_1e^{-C(h)}\leq \alpha_2+\frac{R}{\mu}\alpha_2e^{C(h)}\int_s^t\frac{1}{Y_s(r)}dr.
\end{align*}
When $0\leq t-s\leq h$, applying Gronwall's inequality, we have
\begin{align}\label{Y-s}
  &\,\frac{1}{Y_s(t)}\leq \frac{\alpha_2 C(h)}{\alpha_1}\exp\left\{\frac{R\alpha_2C(h)^2}{\mu\alpha_1}h\right\}\\
  \Longrightarrow&\, Y_s(t)M_s(t)\geq \frac{\alpha_1}{\alpha_2 C(h)}\exp\left\{-\frac{R\alpha_2C(h)^2}{\mu\alpha_1}h-C(h)\right\}:=\tilde{C}(h)>0.
\end{align}

To derive the lower bound on $v(t,x)$, we consider two case:
\begin{itemize}
  \item When $0\leq t\leq \min\{h,T\}$, we take $s=0$, then applying \eqref{M-s} and \eqref{Y-s} yields
  \begin{align}\label{v-lower-1}
    v(t,x)\geq Y_0(t)M_0(t,x)v_0(x)\geq  \tilde{C}(h)\underline{V}.
  \end{align}
  \item When $h\leq t\leq T$, we take $s=t-h$, then applying \eqref{M-s} and \eqref{Y-s} yields
  \begin{align}\label{v-lower-2}
    v(t,x)\geq \frac{R}{\mu}\int_{t-h}^tY_{r}(t)M_{r}(t,x)\theta(r,x)dr\geq \frac{R}{\mu}\tilde{C}(h)\int_{t-h}^t\theta(r,x)dr.
  \end{align}
  Using \eqref{theta-lower} and Jensen's inequality to the convex function $e^x$, we have
  \begin{align}\label{v-lower-3}
    \int_{t-h}^t\theta(r,x)dr&=\int_{t-h}^t\exp(\ln\theta(r,x))dr\notag\\
    &\geq \exp\left\{\int_{t-h}^t\ln\theta(r,x)dr\right\}\notag\\
    &\geq \exp\left\{\int_{t-h}^t\ln\alpha_1dr-\sqrt{\frac{\alpha_2}{\kappa}}\left(\int_{t-h}^t\int_{\bbR}\frac{\kappa\theta_x^2}{v\theta^2}dxdr\right)^{\frac12}\right\}\notag\\
    &\geq \alpha_1h\exp\left\{-\sqrt{\frac{\alpha_2\max\{e_2,e_3\}}{h\kappa}}\right\}
  \end{align}
\end{itemize}
Combining \eqref{v-lower-1} and \eqref{v-lower-3}, we obtain $\underline{v}_h$ which depends only on the initial data and $h$, and is independent of $T$.
In particular, \eqref{v-lower-2} does not involve the pointwise lower bound of
$v(t-h,\cdot)$, so the same $\underline v_h$ can be used at each step
of the continuation. Together with \eqref{upper1}, this completes the proof of the lemma.
\end{proof}

\begin{remark}\label{rem-fixed-time-length}
It is worth noting that $h>0$ can be chosen arbitrarily and then fixed. The lower bound $\underline{v}_h$ depends on $h$, the initial energy and
the initial lower bound $\underline V$, with the physical parameters fixed. In particular, taking $h=1$, we obtain a positive lower bound
independent of time. This choice does not require the local existence time to be one.
\end{remark}
\subsection{Pointwise estimate of temperature}
Now, we derive the uniform upper bound on $\theta(t,x)$.
\begin{lemma}
  \label{lem_the1}
There exists positive constant $C$ independent of $T$ such that
  \begin{equation}\label{lem2.4.1}
    \sup_{0\leq t\leq T}\int_{\mathbb{R}}\left[(\theta-1)^2+u^4+v_x^2\right]
    +\int_0^T\int_{\mathbb{R}}\left[(1+\theta+u^2)u_{x}^2+\theta v_x^2+\chi\Phi_x^2+(1-ve^{-\chi\Phi})^2\right]
    \leq C.
  \end{equation}
\end{lemma}

\begin{proof}
  We divide the proof into six steps.

  \noindent{\bf Step 1}.
Firstly, for each $t\in[0,T]$ and $a>0$, we denote
  $$
  \Omega_a(t):=\{x\in\mathbb{R}: \theta(t,x)>a\}.
  $$
Then it follows from \eqref{vol_0} and \eqref{boun_v} that
\begin{equation}\label{pro2.4.1}
  \sup_{0\leq t\leq T}\int_{\bbR}(1-v)^2+\sup_{0\leq t\leq T}\int_{\mathbb{R}\backslash\Omega_a(t)}(\theta-1)^2\leq C\sup_{0\leq t\leq T}\int_{\mathbb{R}}\eta\leq C.
\end{equation}

  \noindent{\bf Step 2}.
We now estimate the integral $\int_0^T\int_{\Omega_a(t)}\theta_x^2$.
For this purpose, we multiply $\eqref{nsp}_3$ and $\eqref{nsp}_2$
by $(\theta-2)_+:=\max\{\theta-2,0\}$ and $2u(\theta-2)_+$, respectively, then integrate the resultant
over $(0,t)\times \mathbb{R}$ to obtain
\begin{align}
      &\int_{\mathbb{R}}\left[\frac{c_v}{2}(\theta-2)_+^2
      +u^2(\theta-2)_+\right]
      +\int_0^t\int_{\Omega_2(s)}\frac{\kappa\theta_x^2}{v}
      +\int_0^t\int_{\mathbb{R}}\frac{\mu u_{x}^2}{v}(\theta-2)_+\nonumber\\
      =&\int_{\mathbb{R}}\left[\frac{c_v}{2}(\theta_0-2)_+^2
      +u_{0}^2(\theta_0-2)_+\right]+\underbrace{R\int_0^t\int_{\mathbb{R}}\frac{\theta}{v}u_x(\theta-2)_+}_{\mathcal{J}_1}
      +\underbrace{2R\int_0^t\int_{\Omega_2(s)}\frac{\theta}{v}u\theta_x}_{\mathcal{J}_2}\nonumber\\
      &\underbrace{-2\int_0^t\int_{\Omega_2(s)}\frac{\mu u_x}{v}u\theta_x}_{\mathcal{J}_3}
      +\underbrace{\frac{1}{c_v}\int_0^t\int_{\Omega_2(s)}u^2\left(\frac{\mu u_x^2}{v}-\frac{R\theta}{v}u_x\right)}_{\mathcal{J}_4}\nonumber\\
      &+\underbrace{\frac{\kappa}{c_v}\int_0^t\int_{\Omega_2(s)}u^2\left(\frac{\theta_x}{v}\right)_x}_{\mathcal{J}_5}
      +\underbrace{\int_0^t\int_{\bbR}2\frac{\Phi_x}{v}u(\theta-2)_+}_{\mathcal{J}_6}.
      \label{pro2.4.4}
\end{align}
In fact, $\mathcal{J}_1$--$\mathcal{J}_5$ can be estimated similarly as Lemma 2.2 in \cite{MR3466845}. Here, we just list the results.
\begin{align}\label{J1-J5}
  \sum_{i=1}^5|\mathcal{J}_i|\leq~ \eta\int_0^t\int_{\mathbb{R}}\left(u_{x}^2(\theta-2)_++\theta_x^2\right)
  +C(\eta)\int_0^t\sup_{x\in\mathbb{R}}(\theta-\frac32)_+^2+C(\eta)\int_0^t\int_{\mathbb{R}}u^2u_x^2,
\end{align}
where $\eta>0$ is a suitably small constant.
For the last term $\mathcal{J}_6$, we will use different method to control it dependent on the values of $\chi$.

When $\chi=0$, we deduce from Lemma \ref{lem-2.2} that
\begin{align}\label{J-6}
  \mathcal{J}_6=&\int_0^t\int_{\bbR}2\frac{\Phi_x}{v}u(\theta-2)_+
  =-\int_0^t\int_{\bbR}\left(\left(\frac{\Phi_x}{v}\right)^2\right)_t(\theta-2)_+\nonumber\\
  =&-\int_{\bbR}\left(\frac{\Phi_x}{v}\right)^2(\theta-2)_++\int_{\bbR}\left(\frac{\Phi_{0x}}{v_0}\right)^2(\theta_0-2)_++\int_0^t\int_{\Omega_2(s)}\left(\frac{\Phi_x}{v}\right)^2\theta_t.
\end{align}
To estimate the last term on the right-hand side of \eqref{J-6}, we substitute $\theta_t$ defined in $\eqref{nsp}_3$ to derive
\begin{equation*}
  \int_0^t\int_{\Omega_2(s)}\left(\frac{\Phi_x}{v}\right)^2\theta_t=\frac{1}{c_v}\int_0^t\int_{\Omega_2(s)}\left(\frac{\Phi_x}{v}\right)^2\left[\frac{\mu u_x^2}{v}-\frac{R\theta}{v}u_x+\left(\frac{\kappa\theta_x}{v}\right)_x\right].
\end{equation*}
Using H\"{o}lder's inequality, we have
\begin{align}
  \left|\int_0^t\int_{\Omega_2(s)}\left(\frac{\Phi_x}{v}\right)^2\frac{\mu u_x^2}{v}\right|\leq~ & C\int_0^t\left\|\frac{\Phi_x}{v}\right\|_{L^{\infty}}^2\|u_x\|_{L^2}^2\nonumber\\
  \leq~& C\int_0^t\left\|\frac{\Phi_x}{v}\right\|_{L^{2}}\|1-v\|_{L^2}\|u_x\|_{L^2}^2\nonumber\\
  \leq~& \eta\int_0^t\int_{\mathbb{R}}u_x^2\theta+C(\eta)\int_0^t\int_{\mathbb{R}}\frac{u_x^2}{\theta}.
\end{align}
By Cauchy's inequality and \eqref{vol_0}, we obtain
\begin{align}
  \left|\int_0^t\int_{\Omega_2(s)}\left(\frac{\Phi_x}{v}\right)^2\frac{R\theta}{v}u_x\right|\leq~ &
  C\int_0^t\int_{\Omega_2(s)}\left(\left(\frac{\Phi_x}{v}\right)^2\theta^2+\left(\frac{\Phi_x}{v}\right)^2u_x^2\right)\nonumber\\
  \leq~& C\int_0^t\sup_{x\in\bbR}(\theta-\frac{3}{2})_+^2+\eta\int_0^t\int_{\mathbb{R}}u_x^2\theta+C(\eta)\int_0^t\int_{\mathbb{R}}\frac{u_x^2}{\theta}
\end{align}
and
\begin{align}
  \int_0^t\int_{\Omega_2(s)}\left(\frac{\Phi_x}{v}\right)^2\left(\frac{\kappa\theta_x}{v}\right)_x\leq~&
  \left|2\int_0^t\int_{\Omega_2(s)}\frac{\Phi_x}{v}(1-v)\frac{\kappa\theta_x}{v}\right|\nonumber\\
  \leq~& \eta\int_0^t\int_{\bbR}\theta_x^2+C(\eta)\int_0^t\int_{\bbR}(1-v)^2.\label{pro2.4.5}
\end{align}
With above inequalities, we get from $\eqref{boun_v}$ that
\begin{align}
    &\int_{\mathbb{R}}\left((\theta-2)_+^2+\left(\frac{\Phi_x}{v}\right)^2(\theta-2)_++u^2(\theta-2)_+\right)
      +\int_0^t\int_{\Omega_2(s)}\left(\theta_x^2
      +u_x^2(\theta-2)_+\right)\nonumber\\
    \leq~& C
      +C\int_0^t\int_{\mathbb{R}}u^2u_x^2
      +C\int_0^t\sup_{\mathbb{R}}(\theta-\frac32)_+^2+C\int_0^t\int_{\bbR}(1-v)^2.\label{pro2.4.6}
\end{align}

When $\chi=1$, we deduce from $\eqref{boun_v}$ that
\begin{align}\label{J-6-1}
  |\mathcal{J}_6|\leq ~& \eta\int_0^T\int_{\bbR}\Phi_x^2+C(\eta)\int_0^T\int_{\bbR}u^2(\theta-2)^2_+\notag\\
  \leq ~& \eta\int_0^T\int_{\bbR}\Phi_x^2+C(\eta)\int_{0}^T\sup_{x\in\mathbb{R}}(\theta-\frac32)_+^2.
\end{align}
Combining \eqref{J-6-1} and \eqref{pro2.4.7}, and choosing $\eta$
sufficiently small, we have
\begin{align}
      &\int_{\mathbb{R}}\left((\theta-2)_+^2+u^2(\theta-2)_+\right)
      +\int_0^t\int_{\Omega_2(s)}\left(\theta_x^2
      +u_x^2(\theta-2)_+\right)\nonumber\\
      \leq~& C
      +C\int_0^T\sup_{\mathbb{R}}(\theta-\frac32)_+^2
      +C\int_0^T\int_{\mathbb{R}}u^2u_x^2+\eta\int_0^T\int_{\bbR}\Phi_x^2.\label{lem-tem-4}
\end{align}

\noindent{\bf Step 3}.
We obtain from Corollary \ref{cor_0} that
\begin{align}
 \int_0^t\int_{\mathbb{R}}\left[\theta_x^2+u_x^2\theta\right]
 \leq~&\int_0^t\int_{\Omega_3(s)}
 \left[\theta_x^2+u_x^2(\theta-2)_+\right]
 +\int_0^t\int_{\mathbb{R}\setminus\Omega_3(s)}
 \left[\frac{\theta_x^2}{\theta^2}+\frac{u_x^2}{\theta}\right]\nonumber\\
 \leq~&\int_0^t\int_{\Omega_2(s)}
 \left[\theta_x^2+u_x^2(\theta-2)_+\right]
 +C.\label{pro2.4.7}
\end{align}

Combining \eqref{pro2.4.6}, \eqref{lem-tem-4} and \eqref{pro2.4.7}, and choosing $\eta$
sufficiently small, we have
\begin{align}
      &\int_{\mathbb{R}}\left((\theta-2)_+^2+\left(\frac{\Phi_x}{v}\right)^2(\theta-2)_++u^2(\theta-2)_+\right)
      +\int_0^t\int_{\mathbb{R}}\left(\theta_x^2
      +u_x^2\theta\right)\nonumber\\
      \leq~& C
      +C\int_0^t\sup_{\mathbb{R}}(\theta-\frac32)_+^2
      +C\int_0^t\int_{\bbR}(1-v)^2+C\int_0^t\int_{\mathbb{R}}u^2u_x^2, \quad {\rm for}\ \chi=0. \label{pro2.4.8}
\end{align}
and
\begin{align}
      &\int_{\mathbb{R}}\left((\theta-2)_+^2+u^2(\theta-2)_+\right)
      +\int_0^t\int_{\mathbb{R}}\left(\theta_x^2
      +u_x^2\theta\right)\nonumber\\
      \leq~& C
      +C\int_0^t\sup_{\mathbb{R}}(\theta-\frac32)_+^2
      +\eta\int_0^T\int_{\bbR}\Phi_x^2+C\int_0^t\int_{\mathbb{R}}u^2u_x^2,\quad {\rm for}\ \chi=1. \label{pro2.4.8-1}
\end{align}

\noindent{\bf Step 4}.
To estimate the last term on the right-hand side of \eqref{pro2.4.8} and \eqref{pro2.4.8-1}, we
multiply $\eqref{nsp}_2$ by $u^3$ and then integrate the resulting identity
over $(0,t)\times\mathbb{R}$ to have
\begin{align}
   &\frac{1}{4}\int_{\mathbb{R}}u^4
   +3\mu\int_0^t\int_{\mathbb{R}}\frac{u^2u_x^2}{v}-\frac{1}{4}\int_{\mathbb{R}}u_{0}^4
   \nonumber\\
   =~&
   3R\int_0^t\int_{\mathbb{R}}\frac{u^2u_x(1-v)}{v}
   +3R\int_0^t\int_{\mathbb{R}\backslash\Omega_2(s)}\frac{\theta-1}{v}u^2u_x\nonumber\\
   &+3R\int_0^t\int_{\Omega_2(s)}\frac{\theta-1}{v}u^2u_x
   +\int_0^t\int_{\mathbb{R}}\frac{\Phi_x}{v}u^3.\label{pro2.4.9}
\end{align}
The first three terms on the righthand-side of \eqref{pro2.4.9} can be handled by Cauchy's inequality, we have
\begin{align}
    &\left|\int_0^t\int_{\mathbb{R}}\frac{u^2u_x(1-v)}{v}\right|
   +\left|\int_0^t\int_{\mathbb{R}\backslash\Omega_2(s)}\frac{\theta-1}{v}u^2u_x\right|+\left|\int_0^t\int_{\Omega_2(s)}\frac{\theta-1}{v}u^2u_x\right|\nonumber\\
    \leq~& C(\eta)+\eta\int_0^t\int_{\mathbb{R}}\theta u_x^2+\eta\int_0^t\int_{\mathbb{R}}u^2u_x^2
  +C(\eta)\int_0^t\sup_{x\in\mathbb{R}}\left(\theta-\frac32\right)_+^2.\label{pro2.4.10}
\end{align}
For the last term on the right-hand side of \eqref{pro2.4.9}, we shall handle it in different way as well.

When $\chi=0$,
according to Lemma \ref{lem-2.2}, we obtain
\begin{align*}
  \int_0^t\int_{\mathbb{R}}\frac{\Phi_x}{v}u^3=&-\frac{1}{2}\int_0^t\int_{\mathbb{R}}\partial_t\left(\frac{\Phi_x}{v}\right)^2u^2\\
  =&-\frac{1}{2}\int_{\mathbb{R}}\left(\frac{\Phi_x}{v}\right)^2u^2+\frac{1}{2}\int_{\mathbb{R}}\left(\frac{\Phi_{0x}}{v_0}\right)^2u_0^2
  +\int_0^t\int_{\mathbb{R}}\left(\frac{\Phi_x}{v}\right)^2uu_t
\end{align*}
and
\begin{eqnarray*}
  \int_0^t\int_{\mathbb{R}}\left(\frac{\Phi_x}{v}\right)^2uu_t=\int_0^t\int_{\mathbb{R}}\left(\frac{\Phi_x}{v}\right)^2u\left(-\left(\frac{R\theta}{v}\right)_x+\left(\frac{\mu u_x}{v}\right)_x+\frac{\Phi_x}{v}\right):=\sum_{p=1}^{3} \mathcal{K}_p.
\end{eqnarray*}
For $\mathcal{K}_3$, using Lemma \ref{lem-2.2}, we get
\begin{eqnarray}\label{K-3}
  \mathcal{K}_3=-\frac{1}{2}\int_0^t\int_{\mathbb{R}}\left(\frac{\Phi_x}{v}\right)^2\partial_t\left(\frac{\Phi_x}{v}\right)^2
  =-\frac{1}{4}\int_{\mathbb{R}}\left(\frac{\Phi_x}{v}\right)^4+\frac{1}{4}\int_{\mathbb{R}}\left(\frac{\Phi_{0x}}{v_0}\right)^4.
\end{eqnarray}
For $\mathcal{K}_2$, by integrating by parts, Cauchy's inequality, \eqref{pro2.4.1} and \eqref{vol_0}, we obtain
\begin{align}
  |\mathcal{K}_2|=~&\left|\int_0^t\int_{\mathbb{R}}\frac{\Phi_x}{v}(1-v)u\frac{\mu u_x}{v}+\int_0^t\int_{\mathbb{R}}\left(\frac{\Phi_x}{v}\right)^2u_x\frac{\mu u_x}{v}\right|\nonumber\\
  \leq~& \eta\int_0^t\int_{\mathbb{R}}u^2u_x^2+C(\eta)\int_0^t\int_{\mathbb{R}}\left(\frac{\Phi_x}{v}\right)^2(1-v)^2+\int_0^t\left\|\frac{\Phi_x}{v}\right\|_{L^2}\|(1-v)\|_{L^2}\|u_x\|_{L^2}^2\nonumber\\
  \leq~&
  \eta\int_0^t\int_{\mathbb{R}}u^2u_x^2+C(\eta)\int_0^t\int_{\mathbb{R}}(1-v)^2+\eta\int_0^t\int_{\mathbb{R}}\theta u_x^2+C(\eta).\label{K-2}
\end{align}
For $\mathcal{K}_1$, we shall divide it into three parts as follows
\begin{align}\label{K-1}
  \mathcal{K}_1=~&R\int_0^t\int_{\mathbb{R}}\left(\left(\frac{\Phi_x}{v}\right)^2u\right)_x\frac{1-v}{v}+R\int_0^t\int_{\Omega_2(s)}\left(\left(\frac{\Phi_x}{v}\right)^2u\right)_x\frac{\theta-1}{v}\\
  &+R\int_0^t\int_{\mathbb{R}\backslash\Omega_2(s)}\left(\left(\frac{\Phi_x}{v}\right)^2u\right)_x\frac{\theta-1}{v}:=R\sum_{p=1}^{3} \mathcal{K}_{1p}.
\end{align}
For $\mathcal{K}_{11}$,
\begin{align}
  |\mathcal{K}_{11}|\leq~& \int_0^t\int_{\mathbb{R}}u_x^2+\int_0^t\int_{\mathbb{R}}(1-v)^2+\int_0^t\|1-v\|_{L^2}^{\frac{3}{2}}\|u_x\|_{L^2}^{\frac{1}{2}}\nonumber\\
  \leq~& \eta\int_0^t\int_{\mathbb{R}}\theta u_x^2+\int_0^t\int_{\mathbb{R}}(1-v)^2+C(\eta).\label{K-1-1}
\end{align}
For $\mathcal{K}_{12}$,
\begin{align}
  |\mathcal{K}_{12}|\leq~& \int_0^t\int_{\mathbb{R}}\left(\frac{\Phi_x}{v}\right)^2u_x^2+\int_0^t\int_{\mathbb{R}}\left[\left(\frac{\Phi_x}{v}\right)^2+u^2\right](\theta-1)^2+\int_0^t\int_{\mathbb{R}}\left(\frac{\Phi_x}{v}\right)^2(1-v)^2\nonumber\\
  \leq~& \eta\int_0^t\int_{\mathbb{R}}\theta u_x^2+\int_0^t\int_{\mathbb{R}}(1-v)^2+\int_0^t\sup_{x\in\mathbb{R}}\left(\theta-\frac{3}{2}\right)^2_++C(\eta).\label{K-1-2}
\end{align}
For $\mathcal{K}_{13}$,
\begin{align}
  |\mathcal{K}_{13}|\leq~& \int_0^t\|1-v\|_{L^2}^{\frac{3}{2}}\|u_x\|_{L^2}^{\frac{1}{2}}+\int_0^t\|1-v\|_{L^2}\|u_x\|_{L^2}\nonumber\\
  \leq~& \eta\int_0^t\int_{\mathbb{R}}\theta u_x^2+\int_0^t\int_{\mathbb{R}}(1-v)^2+C(\eta).\label{K-1-3}
\end{align}

When $\chi=1$, according to \eqref{boun_v}, we have
\begin{align}\label{lem-tem-6}
  \left|\int_0^T\int_{\mathbb{R}}\frac{\Phi_x}{v}u^3\right|\leq~& C\int_0^T\|\Phi_x\|_{L^2}\|u\|^2_{L^2}\|u_x\|_{L^2}\nonumber\\
  \leq~& \eta\int_0^T\|\Phi_x\|_{L^2}^2+C(\eta)\int_0^T\|u_x\|_{L^2}^2\nonumber\\
  \leq~& \eta\int_0^T\|\Phi_x\|_{L^2}^2+\eta\int_0^T\int_{\bbR}\theta u_x^2+C(\eta)
\end{align}

Plugging \eqref{pro2.4.10}--\eqref{lem-tem-6} into \eqref{pro2.4.9}, combining with \eqref{pro2.4.8} and \eqref{pro2.4.8-1}, and taking $\eta$
sufficiently small, we derive
\begin{align}\label{pro2.4.16}
    &\int_{\mathbb{R}}
    \left[\left(\theta-2\right)_+^2+\left(\frac{\Phi_x}{v}\right)^2(\theta-2)_++u^4+\left(\frac{\Phi_x}{v}\right)^2u^2\right]
    +\int_0^t\int_{\mathbb{R}}\left[\theta_x^2+u_x^2\left(1+\theta+u^2\right)\right]\nonumber\\
    \leq~& C+C\int_0^t\sup_{\mathbb{R}}\left(\theta-\frac32\right)_+^2+C\int_0^t\int_{\mathbb{R}}(1-v)^2,\quad {\rm for}\ \chi=0.
 \end{align}
and
\begin{align}\label{lem-tem-8}
    &\int_{\mathbb{R}}
    \left[\left(\theta-2\right)_+^2+u^4\right]
    +\int_0^T\int_{\mathbb{R}}\left[\theta_x^2+u_x^2\left(1+\theta+u^2\right)\right]\nonumber\\
    \leq~& C+C\int_0^T\sup_{\mathbb{R}}\left(\theta-\frac32\right)_+^2+\eta\int_0^T\|\Phi_x\|_{L^2}^2,\quad {\rm for}\ \chi=1.
 \end{align}

\noindent{\bf Step 5}.
For the second term on the right-hand side of \eqref{pro2.4.16} and \eqref{lem-tem-8},
according to the fundamental theorem of calculus, we have from Corollary \ref{cor_0} that
\begin{align}
    \int_0^T\sup_{\mathbb{R}}\left(\theta-\frac32\right)_+^2
    \leq \eta\int_0^T\int_{\mathbb{R}}\theta_x^2+C(\eta).\label{pro2.4.17}
\end{align}

\noindent{\bf Step 6}.
It remains to estimate the last term on the right-hand side of \eqref{pro2.4.16} and \eqref{lem-tem-8}.
Multiplying \eqref{nsp}$_1$ by $\frac{v_x}{v}$ and integrating over $\bbR$, we have
\begin{align*}
  &\frac{\mu}{2}\int_{\bbR}\left(\frac{v_x}{v}\right)^2+\int_0^T\int_{\bbR}\frac{R\theta v_x^2}{v^3}-\frac{\mu}{2}\int_{\bbR}\left(\frac{v_{0x}}{v_0}\right)^2\\
  =&\int_0^T\int_{\bbR}\frac{R\theta_xv_x}{v^2}+\int_0^T\int_{\bbR}\frac{u_tv_x}{v}-\int_0^T\int_{\bbR}\frac{\Phi_xv_x}{v^2}
  :=\sum_{p=1}^3\mathcal{L}_{p}.
\end{align*}
Here, $\mathcal{L}_{i}(i=1,2)$ can be estimated by Cauchy's inequality as follows:
\begin{align}
  |\mathcal{L}_{1}|+|\mathcal{L}_{2}|\leq \frac{1}{4}\int_0^T\int_{\bbR}\frac{R\theta v_x^2}{v^3}+\eta\int_0^T\int_{\bbR}\theta_x^2+\frac{\mu}{4}\int_{\bbR}\left(\frac{v_x}{v}\right)^2+\eta\int_0^T\int_{\bbR}\theta u_x^2+C(\eta).\label{L-1}
\end{align}

$\mathcal{L}_{3}$ is crucial in our analysis. Actually, using integration by parts yields
\begin{align}\label{L-3-1}
  \mathcal{L}_{3}&=-\int_0^T\int_{\bbR}\frac{\Phi_x}{v}(\ln v+\ln(e^{-\chi\Phi}))_x-\int_0^T\int_{\bbR}\chi\frac{\Phi_x^2}{v}\notag\\
  &=\int_0^T\int_{\bbR}(1-ve^{-\chi\Phi})\ln(ve^{-\chi\Phi})-\int_0^T\int_{\bbR}\chi\frac{\Phi_x^2}{v},
\end{align}
where
\begin{align*}
  \int_0^T\int_{\bbR}(1-ve^{-\chi\Phi})\ln(ve^{-\chi\Phi})\leq 0.
\end{align*}

Combining \eqref{L-1}$-$\eqref{L-3-1}, it depends on the uniform bounds of $v(t,x)$ that
\begin{eqnarray}\label{v-x}
  \int_{\bbR}v_x^2+\int_0^T\int_{\bbR}\theta v_x^2+\int_0^T\int_{\bbR}(1-ve^{-\chi\Phi})^2+\int_0^T\int_{\bbR}\chi\Phi_x^2
  \leq \eta\int_0^T\int_{\bbR}\theta_x^2+\eta\int_0^T\int_{\bbR}\theta u_x^2+C(\eta).
\end{eqnarray}
Plug \eqref{pro2.4.17}, \eqref{v-x} into \eqref{pro2.4.16} and \eqref{lem-tem-8}, and choose $\eta>0$
suitably small to obtain \eqref{lem2.4.1}.
This completes the proof of the lemma.
\end{proof}

We obtain the upper bound for the temperature uniformly in both time and
space in the next lemma,
by combining Lemma \ref{lem_the1} and some desired uniform estimates
on the spatial derivatives of $(u,\theta)$.
\begin{lemma}\label{lem_the2}
There exists a positive constant $C$ independent of $T$ such that
  \begin{gather}\label{lem2.4.3}
   \theta(t,x)\leq \overline{\theta} \quad {\rm for\ all}\ (t,x)\in[0,T]\times\mathbb{R},\\
    \sup_{0\leq t\leq T}\|(v-1,u,\theta-1,\chi\Phi)(t)\|_{1}^2
    +\int_0^T[\|(\theta v_x,\chi\Phi_x)\|_{L^2}^2+\|(u_x,\theta_{x})(t)\|_1^2]dt
    \leq C.\label{lem2.3.2}
  \end{gather}
\end{lemma}

\begin{proof}
Multiplying $\eqref{nsp}_2$ by $-u_{xx}$ and integrating over $(0,T)\times\mathbb{R}$, we have
\begin{eqnarray*}
  \frac12\int_{\mathbb{R}}u_x^2+\mu\int_0^T\int_{\mathbb{R}}\frac{u_{xx}^2}{v}=\frac12\int_{\mathbb{R}}u_{0x}^2+
  \int_0^T\int_{\mathbb{R}}R\left(\frac{\theta}{v}\right)_xu_{xx}+\int_0^T\int_{\mathbb{R}}\mu\frac{u_x}{v^2}v_xu_{xx}-\int_0^T\int_{\mathbb{R}}\frac{\Phi_x}{v}u_{xx}.
\end{eqnarray*}
Next, we only show the estimates for terms contain $\Phi_x$ for brevity, since the rest proof can be followed in the same way as Lemma 2.3 in \cite{MR3466845}.
Integrating by parts and applying Lemma \ref{lem-2.2} imply that
\begin{align}
  -\int_0^T\int_{\mathbb{R}}\frac{\Phi_x}{v}u_{xx}=&\int_0^T\int_{\mathbb{R}}(1-v)u_x\nonumber\\
  \leq& -\frac12\int_{\mathbb{R}}(1-v)^2+\frac12\int_{\mathbb{R}}(1-v_0)^2,\quad {\rm for}\ \chi=0,\label{u-x-2}
\end{align}
and
\begin{align}
  \int_0^T\int_{\mathbb{R}}\frac{\Phi_x}{v}u_{xx}
  \leq \eta\int_0^T\int_{\mathbb{R}}u_{xx}^2+C(\eta)\int_0^T\int_{\mathbb{R}}\Phi_x^2,\quad {\rm for}\ \chi=1.\label{u-x-2-1}
\end{align}
This completes the proof of the lemma.
\end{proof}

Finally, to derive the existence of global solution by using continuity argument, we give a positive lower bound of $\theta(t,x)$ on $[0,T]$ according to maximum principle.
\begin{lemma}[c.f.\cite{MR3564590}]\label{lem-theta-low}
  There exists a positive constant $C$ independent of $T$ such that
  \begin{equation*}
    \inf_{\mathbb{R}}\theta(t,\cdot)\geq \frac{\inf_{\mathbb{R}}\theta(s,\cdot)}
	{C\inf_{\mathbb{R}}\theta(s,\cdot)(t-s)+1}
	\quad
	{\rm for\ all}\
	0\leq s\leq t\leq T.
  \end{equation*}
\end{lemma}

\subsection{Proof of Theorem \ref{main-theorem}}
The global existence can be derived by the local existence in Lemma \ref{lem-2.1}, and the a priori estimates in Lemmas \ref{lem_the1}-\ref{lem-theta-low} combining
with the continuity argument below.

Fix $h>0$ and denote by $T_{\max}$ the maximal existence time of the local solution. For every
$T<T_{\max}$, Lemma \ref{lem_boun} yields
\begin{equation}\label{continuation-v}
 \underline v_h\leq v(t,x)\leq\overline v,
 \qquad (t,x)\in[0,T]\times\mathbb R,
\end{equation}
where $\underline v_h$ is independent of $T$. Moreover, it follows from Lemmas \ref{lem_the1} and
\ref{lem_the2} that
\begin{equation}\label{continuation-H1}
 \sup_{0\leq t\leq T}\|(v-1,u,\theta-1,\chi\Phi)(t)\|_1\leq C(h),
 \qquad \theta(t,x)\leq\overline\theta.
\end{equation}
Here $C(h)$ depends only on the initial data and $h$, and is independent of $T$.

For the electric potential, Corollary \ref{cor_0} gives $\left\|\frac{\Phi_x}{v}\right\|_{L^2}\leq C(h)$.
When $\chi=0$, we have $\left\|\left(\frac{\Phi_x}{v}\right)_x\right\|_{L^2}=\|1-v\|_{L^2}\leq C(h)$. When $\chi=1$, applying Lemma \ref{lem2.3} and \eqref{continuation-H1}, we obtain
\begin{equation*}
 \left\|\left(\frac{\Phi_x}{v}\right)_x\right\|_{L^2}\leq\|1-v\|_{L^2}+\overline v\|1-e^{-\Phi}\|_{L^2}
 \leq\|1-v\|_{L^2}+C(h)\|\Phi\|_{L^2}\leq C(h).
\end{equation*}
Hence, Sobolev's inequality implies $\left\|\frac{\Phi_x}{v}\right\|_{L^\infty}\leq C(h)$.
Applying $\left(\frac{\Phi_x}{v}\right)_x=1-ve^{-\chi\Phi}$, we get
\begin{equation}\label{continuation-potential}
 \begin{cases}
 \|\Phi_x(t)\|_1\leq C(h),&\chi=0,\\
 \|\Phi(t)\|_2\leq C(h),&\chi=1.
 \end{cases}
\end{equation}

On the other hand, let $\theta_*:=\inf_{\mathbb R}\theta_0>0$. For any finite $S>0$, Lemma \ref{lem-theta-low} yields
\begin{equation}\label{continuation-theta}
 \theta(t,x)\geq m_\theta(S):=
 \frac{\theta_*}{1+C\theta_*S}>0,
 \qquad 0\leq t<\min\{S,T_{\max}\}.
\end{equation}
Combining \eqref{continuation-v}--\eqref{continuation-theta}, we can apply Lemma \ref{lem-2.1} at any time $t_0<\min\{S,T_{\max}\}$ with
a common local existence time $\tau_{h,S}>0$, which is independent of $t_0$. The Poisson equation and the far-field conditions hold at
$t=t_0$ since the data are taken from the existing solution. If $T_{\max}<\infty$, we take $S=T_{\max}$ and choose $t_0$ such that
$T_{\max}-t_0<\tau_{h,S}/2$. The local solution starting at $t_0$ then exists beyond $T_{\max}$ and agrees with the previous solution
on their common interval by uniqueness. This contradicts the definition of $T_{\max}$. Therefore, $T_{\max}=\infty$.

The large time behavior and the uniform lower bound for $\theta(t,x)$ are established in the following lemma which completes the proof of Theorem \ref{main-theorem}.
\begin{lemma}\label{lem-2.9}
For any $p>2$, it holds that
\begin{itemize}
  \item When $\chi=0$,
  \begin{equation*}
  \lim_{t\rightarrow\infty}\left(\|(v-1,u,\theta-1,\Phi_x)\|_{L^p}+\|(1-v,v_x,u_x,\theta_x)\|_{L^2}\right)=0.
\end{equation*}
  \item When $\chi=1$,
  \begin{equation*}
  \lim_{t\rightarrow\infty}\left(\|(v-1,u,\theta-1,\Phi,\Phi_x)\|_{L^p}+\|(v_x,u_x,\theta_x,\Phi_x)\|_{L^2}\right)=0.
\end{equation*}
\end{itemize}
Moreover, there exists a positive constant $C_8$ such that
\begin{equation*}
  C_8^{-1}\leq \theta\leq C_8.
\end{equation*}
\end{lemma}

\begin{proof}
Combining Lemma \ref{lem_the1} and Lemma \ref{lem_the2}, we derive that
\begin{align*}
  &\int_0^{\infty}\left(\|u_x(t,x)\|_{L^2}^2+\left|\frac{d}{dt}\|u_x(t,x)\|_{L^2}^2\right|\right)dt\\
  &+\int_0^{\infty}\left(\|\theta_x(t,x)\|_{L^2}^2+\left|\frac{d}{dt}\|\theta_x(t,x)\|_{L^2}^2\right|\right)dt\leq C,
\end{align*}
which gives
\begin{equation*}
  \lim_{t\rightarrow\infty}\|(u_x,\theta_x)\|_{L^2}=0,\quad  \lim_{t\rightarrow\infty}\|(u,\theta-1)\|_{L^{\infty}}=0.
\end{equation*}
This concludes that there exists some $T_0>0$ such that for all $x\in[T_0,\infty)\times\bbR$,
\begin{equation*}
  \frac12\leq\theta\leq \frac32.
\end{equation*}
Hence, it derives by Lemma \ref{lem_the2} that
\begin{equation*}
  \int_{0}^{\infty}\left(\|v_x(t,x)\|_{L^2}^2+\left|\frac{d}{dt}\|v_x(t,x)\|_{L^2}^2\right|\right)dt\leq C,
\end{equation*}
which implies
\begin{equation*}
  \lim_{t\rightarrow\infty}\|v_x\|_{L^2}=0,\quad  \lim_{t\rightarrow\infty}\|v-1\|_{L^{\infty}}=0.
\end{equation*}

On the other hand, since we can not control $\int_{0}^{\infty}\int_{\bbR}\Phi_x^2(t,x)$ when $\chi=0$, it is hopeless to obtain the large time behavior for $\|\Phi_x(t,x)\|_{L^2}$.
Instead, we could derive the asymptotic behavior for $\Phi_{xx}(t,x)$. It depends on equation $\eqref{nsp}_4$ and \eqref{lem2.4.1} that
\begin{equation*}
  \int_{0}^{\infty}\left(\left\|\left(\frac{\Phi_{x}}{v}\right)_x\right\|_{L^2}^2+\left|\frac{d}{dt}\left\|\left(\frac{\Phi_{x}}{v}\right)_x\right\|_{L^2}^2\right|\right)dt\leq C,
\end{equation*}
which yields
\begin{equation*}
  \lim_{t\rightarrow\infty}\left\|\left(\frac{\Phi_{x}}{v}\right)_x\right\|_{L^2}=\lim_{t\rightarrow\infty}\|v-1\|_{L^{2}}=0,\quad  \lim_{t\rightarrow\infty}\left\|\Phi_{x}\right\|_{L^{\infty}}=0.
\end{equation*}

When $\chi=1$, due to the damping effect induced by Boltzmann relation, we
apply $\partial_t$ to \eqref{nsp}$_{4}$, multiply the resultant by $\Phi_t$ and using H\"{o}lder's inequality as well as Young's inequality to derive (cf. \cite{MR4458894})
\begin{equation*}
  \|(\Phi_t,\Phi_{xt})\|_{L^2}^2\leq \|(v_x,u_x)\|_{L^2}^2,
\end{equation*}
which combines with Lemma \ref{lem_the2} yields
\begin{equation*}
  \int_{0}^{\infty}\left(\|\Phi_x(t,x)\|_{L^2}^2+\left|\frac{d}{dt}\|\Phi_x(t,x)\|_{L^2}^2\right|\right)dt\leq C,
\end{equation*}
so we obtain
\begin{equation}\label{phi-x-1}
  \lim_{t\rightarrow\infty}\|\Phi_x\|_{L^{2}}=0,\quad \lim_{t\rightarrow\infty}\|\Phi_x\|_{L^{\infty}}=0.
\end{equation}
Moreover, we derive from \eqref{cor-lag-2} and Lemma \ref{lem2.3} that
\begin{align}\label{phi-l3}
  \|\Phi\|_{L^2}^2\leq \left\|\frac{\Phi^2}{1-\Phi e^{-\Phi}-e^{-\Phi}}\right\|_{L^{\infty}}\left\|1-\Phi e^{-\Phi}-e^{-\Phi}\right\|_{L^1}\leq C,
\end{align}
where we have used the fact that
\begin{align*}
  \lim_{\Phi\rightarrow0}\frac{\Phi^2}{1-\Phi e^{-\Phi}-e^{-\Phi}}=2.
\end{align*}
Combining \eqref{phi-l3} with \eqref{phi-x-1} and using Sobolev's inequality, we obtain
\begin{align*}
  \lim_{t\rightarrow\infty}\|\Phi\|_{L^{\infty}}\leq C\lim_{t\rightarrow\infty}\|\Phi\|_{L^2}^{\frac12}\|\Phi_x\|_{L^2}^{\frac12}=0.
\end{align*}

Finally, we can derive a lower bound for the temperature in the same way as \cite{MR0468593} for some constant $C_9>0$ as follows,
\begin{equation*}
  C_9^{-1}\exp^{-C_9t}\leq \theta(t,x).
\end{equation*}
By choosing $C_8:=\max\{\overline{\theta},C_9\exp^{C_9T_0}\}$, we finish the proof.
\end{proof}

\section{The Proof of Theorem \ref{main-theorem-1}}
In this section, the proof of Theorem \ref{main-theorem-1} will be presented. Firstly, the local existence of the solutions can be proved similarly as Lemma \ref{lem-2.1}.
Then, we define a functional space of the solutions, for constants $N$, $m_i$, $s$, $t$ ($i=1,2$, $t\geq s$) and $\forall(\tau,x)\in[s,t ]\times\mathbb{R}$,
in the following way.
\begin{equation*}
X(s,t ;m_1,m_2,N):=\left\{(v, u,\theta)\left|
\begin{array}{lll}
&(v-1,u,\theta-1)\in C([s,t];H^3),\\
&\Phi_x\in L^{\infty}([s,t];H^3),\quad {\rm for\ } \chi=0,\\
&\Phi\in L^{\infty}([s,t];H^4),\quad {\rm for\ } \chi=1,\\
&v_x\in L^2(s,t;H^2),\
(u_x,\theta_x)\in L^2(s,t;H^3),\\
&
 \mathcal{E}_{\chi+1}(s,t)\leq N^2,\ v(\tau,x)\geq m_1,\
\theta(\tau,x)\geq m_2,
\end{array}
\right.
\right\}
\end{equation*}
where
\begin{align*}
\mathcal{E}_1(s,t):=\sup_{\tau\in[s,t]}\|(v-1,u,\theta-1)(\tau)\|_3^2+\sup_{\tau\in[s,t]}\|\Phi_x\|_3^2
+\int_s^t\left[\|v_x(\tau)\|_2^2+\|(u_x,\theta_x)(\tau)\|_3^2\right]
\mathrm{d}\tau,\\
\mathcal{E}_2(s,t):=\sup_{\tau\in[s,t]}\|(v-1,u,\theta-1)(\tau)\|_3^2+\sup_{\tau\in[s,t]}\|\Phi\|_4^2
+\int_s^t\left[\|v_x(\tau)\|_2^2+\|(u_x,\theta_x,\Phi_x)(\tau)\|_3^2\right]
\mathrm{d}\tau.
\end{align*}

Next, we will derive certain a priori estimates for the solutions $(v(t,x),u(t,x),\theta(t,x))$ which belongs to $X(0,T;m_1,m_2,N)$. Here $0<m_{i}\leq 1\leq N<+\infty$.
It follows from the Sobolev's inequality that
\begin{align}\label{priori-v}
  m_1\leq v(t,x)\leq 4 N, \quad
m_2\leq \theta(t,x)\leq 4 N\quad
{\rm for\ all}\ (t,x)\in[0,T]\times\mathbb{R}.
\end{align}

\subsection{Pointwise estimate of specific volume}
In this part, we will derive the upper and lower bounds for the specific volume $v(t,x)$ which depends on the temperature $\theta(t,x)$. Firstly, we show an elementary energy estimates as Corollary \ref{cor_0} as follows.
\begin{corollary}\label{cor-1}
Under the assumptions in Theorem \ref{main-theorem-1}, the following energy estimates hold,
\begin{itemize}
  \item When $\chi=0$, we have
  \begin{align}\notag
  &\sup_{0\leq t\leq T}\int_{\bbR}\left[\eta(v,u,\theta)+\frac{1}{2}\left(\frac{\Phi_x}{v}\right)^2\right](t,x)dx
    +\int_0^{T}\int_{\bbR}\left[\frac{\mu(v,\theta) u_{x}^2}{v \theta}
  +\frac{\kappa(v,\theta)\theta_x^2}{v\theta^2}\right]dxds \nonumber\\
  \leq~ &\int_{\bbR}\left[\eta(v,u,\theta)+\frac{1}{2}\left(\frac{\Phi_x}{v}\right)^2\right](0,x)dx.\label{basic-1}
\end{align}
  \item When $\chi=1$, using Lemma \ref{lem2.2}, we obtain
\begin{align}
  &\sup_{0\leq t\leq T}\int_{\bbR}\left[\eta(v,u,\theta)+\frac12\left(\frac{\Phi_x}{v}\right)^2e^{\Phi}+\frac16\frac{\Phi_x^4}{v^3}e^{\Phi}+v(1-\Phi e^{-\Phi}-e^{-\Phi})\right](t,x)dx \nonumber\\
    &\hspace{1cm}+\int_0^{T}\int_{\bbR}\left[\frac{\mu(v,\theta) u_{x}^2}{v \theta}
  +\frac{\kappa(v,\theta)\theta_x^2}{v\theta^2}\right]dxds \nonumber\\
  \leq~ &C(\underline{V},\overline{V})\left[\int_{\bbR}\eta(v,u,\theta)(0,x)dx+\|v_0-1\|_{L^2}^2+\|v_0-1\|_{L^2}^6\right].\label{basic-2}
\end{align}
\end{itemize}
\end{corollary}

Next, we present an estimate on $\left\|\frac{\mu v_x}{v}\right\|_{L^2}$ which is crucial to derive the bounds from above to below of the specific of volume $v(t,x)$.
\begin{lemma}\label{lem-3.1}
Under the assumptions in Theorem \ref{main-theorem-1}, there exists a constant $\varepsilon_1>0$ such that if
\begin{equation}\label{3.1-assump}
  m_2^{-|\alpha|}\leq 2,\quad N^{|\alpha|}\leq 2,\quad \Xi(m_1,m_2,N)|\alpha|\leq \varepsilon_1,
\end{equation}
where $$
	\Xi(m_1,m_2,N):=
	\left[m_1^{-1}+m_2^{-1}+N+
	\sup_{m_1\leq v \leq 4N}h(v)+\sup_{m_1\leq v \leq 4N}h^{-1}(v)+1\right]^{80},$$
then
\begin{align}\label{vx-v}
  \sup_{t\in[0,T]}\left\|\frac{\mu v_x}{v}(t)\right\|_{L^2}^2
  +\int_0^T\int_{\mathbb{R}}\left(\frac{\mu \theta v_x^2}{v^3}+\chi\frac{\mu(e^{\chi\Phi},\theta)}{v}\Phi^2_x+\left(g(v,\theta)-g(e^{\chi\Phi},\theta)\right)\left(ve^{-\chi\Phi}-1\right)\right)
  \lesssim 1+\|\theta\|_{L^{\infty}_{t,x}},
\end{align}
where
\begin{equation*}
  g(v,\theta)=\int_1^v\frac{\mu(z,\theta)}{z}dz\thicksim\theta^{\alpha}\left[\frac{1}{l_1}(v^{l_1}-1)+\frac{1}{l_2}(1-v^{-l_2})\right].
\end{equation*}
\end{lemma}
\begin{proof}
Direct calculations show that
\begin{equation}\label{3.1-1}
	\left(\frac{\mu v_x}{v}\right)_t
	=u_t+\left(\frac{R\theta}{v}\right)_x-\frac{\Phi_x}{v}
	+\frac{\mu_{\theta}}{v}(v_x \theta_t-\theta_x u_x).
\end{equation}
Multiply \eqref{3.1-1} by ${\mu v_x}/{v}$ and integrate with respect to time and space, we have
\begin{equation}\label{3.1-2}
	\begin{aligned}
	&\frac{1}{2}\left\|\frac{\mu v_x}{v}(t)\right\|_{L^2}^2
	+\int_0^t\int_{\mathbb{R}}\frac{\mu\theta v_x^2}{v^3}
	= \frac{1}{2}\left\|\frac{\mu v_{x}}{v}(0)\right\|_{L^2}^2+\int_0^t\int_{\mathbb{R}}\left(\frac{\mu v_x}{v}u\right)_t
	+\int_0^t\int_{\mathbb{R}}\frac{\mu u_x^2}{v} \\
	 &\qquad
	+\int_0^t\int_{\mathbb{R}}\frac{\mu v_x\theta_x}{v^2}
	-\int_0^t\int_{\mathbb{R}}\frac{\mu v_x}{v}\frac{\Phi_x}{v}
    +\int_0^t\int_{\mathbb{R}}
	\frac{\mu_{\theta}}{v^2}(\mu v_x-uv)(v_x \theta_t-\theta_x u_x).
	\end{aligned}
\end{equation}
We deduce from Cauchy's inequality, \eqref{basic-1} and \eqref{transport} that
\begin{equation}\label{3.1-3}
  \left|\int_0^t\int_{\mathbb{R}}\frac{\mu u_x^2}{v}\right|+\left|\int_0^t\int_{\mathbb{R}}\frac{\mu v_x\theta_x}{v^2}\right|\leq \frac{1}{2}\int_0^t\int_{\mathbb{R}}\frac{\mu\theta v_x^2}{v^3}+C\|\theta\|_{L^{\infty}_{t,x}}.
\end{equation}
Then by integrating by parts, we have
\begin{align}\label{3.1-4}
  &\int_0^t\int_{\mathbb{R}}\frac{\mu(v,\theta) v_x}{v}\frac{\Phi_x}{v}
  -\int_0^t\int_{\mathbb{R}}\chi\frac{\mu(e^{\chi\Phi},\theta)}{v}\Phi^2_x+\int_0^t\int_{\mathbb{R}}(g(v,\theta)-g(e^{\chi\Phi},\theta))\frac{\alpha\theta_x\Phi_x}{v\theta}\nonumber\\
  =&\int_0^t\int_{\mathbb{R}}(g(v,\theta)-g(e^{\chi\Phi},\theta))_x\frac{\Phi_x}{v}\nonumber\\
  =&-\int_0^t\int_{\mathbb{R}}(g(v,\theta)-g(e^{\chi\Phi},\theta))(1-ve^{-\chi\Phi})\geq0.
\end{align}
We have used the monotonicity of $g(v,\theta)$ with respect to $v$. We apply Lemma \ref{lem2.3} to obtain
\begin{align*}
  \int_0^t\int_{\mathbb{R}}(g(v,\theta)-g(e^{\chi\Phi},\theta))\frac{\alpha\theta_x\Phi_x}{v\theta}
  \lesssim&\, |\alpha|\chi\int_0^t\int_{\mathbb{R}}\frac{\mu(e^{\chi\Phi},\theta)}{v}\Phi^2_x
  +|\alpha|\Xi(m_1,m_2,N)\int_0^t\int_{\mathbb{R}}\frac{\kappa\theta_x^2}{v\theta^2}\\
  &\,+|\alpha|(1-\chi)\int_0^t\int_{\mathbb{R}}g(v,\theta)(v-1).
\end{align*}

For the last term in the righthand-side of \eqref{3.1-2}, we rewrite it as follows:
\begin{equation}\label{3.1-5}
	(\mu v_x-uv)(v_x \theta_t-\theta_x u_x)
	=uv\theta_x u_x+\mathcal{R}_1v_x+\mathcal{R}_2v_x^2
\end{equation}
with
\begin{align*}
	\mathcal{R}_1&:=-\mu\theta_x u_x
	-\frac{u}{c_v}\left(\kappa_{\theta}\theta_x^2+\kappa\theta_{xx}
	+\mu u_x^2-R\theta u_x\right),\\
	\mathcal{R}_2&:=\mu\theta_t-\frac{u}{c_v}\left[
	\kappa_v\theta_x-\frac{\kappa\theta_x}{v}\right].
\end{align*}
Inserting \eqref{3.1-3}-\eqref{3.1-5} into \eqref{3.1-2}, we obtain
\begin{align} \notag
	&\left\|\frac{\mu v_x}{v}(t)\right\|_{L^2}^2
	+\int_0^t\int_{\mathbb{R}}\left(\frac{\mu\theta v_x^2}{v^3}+\chi\frac{\mu(e^{\chi\Phi},\theta)}{v}\Phi^2_x+\left(g(v,\theta)-g(e^{\chi\Phi},\theta)\right)\left(ve^{-\chi\Phi}-1\right)\right)
	\lesssim 1+\|\theta\|_{L^{\infty}_{t,x}}\\
	\label{3.1-6}
	 &\qquad
	+\left|\int_0^t\int_{\mathbb{R}}\frac{\mu_{\theta}}{v}u\theta_x u_x \right|
	+\left|\int_0^t\int_{\mathbb{R}}\frac{\mu_{\theta}}{v^2}v_x
	\mathcal{R}_1  \right|
	+\left|\int_0^t\int_{\mathbb{R}}\frac{\mu_{\theta}}{v^2}v_x^2
	\mathcal{R}_2 \right|.
\end{align}
The rest estimates can be done in the same way as Lemma 2.2 in \cite{MR3564590}, we just show the results without any detail for brevity.
\begin{align}\label{3.1-7}
	&\left|\int_0^t\int_{\mathbb{R}}\frac{\mu_{\theta}}{v}u\theta_x u_x\right|+\left|\int_0^t\int_{\mathbb{R}}\frac{\mu_{\theta}}
	{v^2}v_x\mathcal{R}_1\right|+\left|\int_0^t\int_{\mathbb{R}}\frac{\mu_{\theta}}{v^2}v_x^2\mathcal{R}_2\right|\nonumber\\
	\lesssim~& |\alpha| N^{\frac32}+\eta\int_0^t\int_{\mathbb{R}}\frac{\mu\theta v_x^2}{v^3}
	+\alpha^2 \Xi(m_1,m_2,N)^{\frac{1}{8}}\left[1
	+\int_0^t\int_{\mathbb{R}}\frac{\kappa\theta_{xx}^2}{v}\right]
    +|\alpha|\Xi(m_1,m_2,N)^{\frac{1}{4}}
	\int_0^t\int_{\mathbb{R}}\frac{\mu\theta v_x^2 }{v^3},
\end{align}
Substituting \eqref{3.1-7} into \eqref{3.1-6} and choosing a suitably small $\varepsilon_1$, we can deduce
\begin{align}\label{3.1-10}
	&\left\|\frac{\mu v_x}{v}(t)\right\|_{L^2}^2
	+\int_0^t\int_{\mathbb{R}}\left(\frac{\mu\theta v_x^2}{v^3}+\chi\frac{\mu(e^{\chi\Phi},\theta)}{v}\Phi^2_x+\left(g(v,\theta)-g(e^{\chi\Phi},\theta)\right)\left(ve^{-\chi\Phi}-1\right)\right)\nonumber\\
	\lesssim& 1+\|\theta\|_{L^{\infty}_{t,x}}
	+|\alpha|\int_0^t\int_{\mathbb{R}}\frac{\kappa\theta_{xx}^2}{v }.
\end{align}

To estimate the last term in the righthand-side of \eqref{3.1-10}, we multiply \eqref{nsp}$_3$ by $\theta_{xx}$, integrate the resultant over $[0,t]\times\mathbb{R}$ and employ Cauchy's inequality, \eqref{3.1-assump} to get
\begin{align} \notag
	\|\theta_x(t)\|_{L^2}^2+\int_0^t\int_{\mathbb{R}}\frac{\kappa\theta_{xx}^2}{v}
	\lesssim&~1+\int_0^t\int_{\mathbb{R}}
	\left[\frac{\theta^2 u_x^2}{v\kappa}+\frac{\kappa_v^2 v_x^2\theta_x^2}{v\kappa}
	+\frac{\kappa\theta_x^2 v_x^2}{v^3}+\frac{\kappa_{\theta}^2\theta_x^4}{v\kappa}
	+\frac{\mu^2 u_x^4}{v\kappa}\right]\\ \notag
	\lesssim&~1+\left\|\frac{v\theta}{\mu}\left(\frac{\theta^2}{v\kappa}
	+\frac{\mu^2u_x^2}{v\kappa}\right)\right\|_{L^{\infty}}\int_0^t\int_{\mathbb{R}}\frac{\mu u_x^2}{v\theta}\\ \label{pro2.2}
	&+\left\|\frac{v\theta^2}{\kappa}\left(\frac{\kappa_v^2v_x^2}{v\kappa}
	+\frac{\kappa v_x^2}{v^3}+\frac{\kappa_{\theta}^2\theta_x^2}{v\kappa}\right)\right\|_{L^{\infty}}
	\int_0^t\int_{\mathbb{R}}\frac{\kappa \theta_x^2}{v\theta^2}\nonumber\\
    \lesssim&~ 1+\Xi(m_1,m_2,N)^{\frac{1}{8}}.
\end{align}

Combining \eqref{3.1-10}, \eqref{pro2.2} and \eqref{3.1-assump},  we obtain
\begin{align} \label{3.1-12}
	\left\|\frac{\mu v_x}{v}(t)\right\|_{L^2}^2
	+\int_0^t\int_{\mathbb{R}}\left(\frac{\mu\theta v_x^2}{v^3}+\chi\frac{\mu(e^{\chi\Phi},\theta)}{v}\Phi^2_x+\left(g(v,\theta)-g(e^{\chi\Phi},\theta)\right)\left(ve^{-\chi\Phi}-1\right)\right)
	\lesssim  1+\|\theta\|_{L^{\infty}_{t,x}}.
\end{align}
This completes the proof of this lemma.
\end{proof}

With Lemma \ref{lem-3.1} in hand, we can use the method introduced by Kanel$'$ (see \cite{MR0227619,MR3225502,MR3564590}) to derive the upper and lower bounds of the specific volume $v(t,x)$ in term of $\|\theta\|_{L^{\infty}_{t,x}}$. We omit the detail for brevity, whose proof can be found in \cite{MR3564590}.
\begin{lemma}\label{lem-3.3}
Under the assumptions in Lemma \ref{lem-3.1}, the specific volume $v(t,x)$ satisfies
 \begin{equation}\label{v-bound}
	\|v\|_{L^{\infty}_{t,x}}\lesssim 1+\|\theta\|_{L^{\infty}_{t,x}}^{\frac1{1+2{\ell_1}}},
	\qquad
	\|v^{-1}\|_{L^{\infty}_{t,x}}\lesssim 1+\|\theta\|_{L^{\infty}_{t,x}}^{\frac1{2{\ell_2}}}.
\end{equation}
\end{lemma}

\subsection{Pointwise estimate of temperature}
In this subsection, we will devote to deriving the upper and lower bounds of the temperature $\theta(t,x)$.
The following $L^2$ estimate of $\theta-1$ is of great importance to accomplish this goal.

\begin{lemma}\label{lem-3.4}
Under the assumptions in Lemma \ref{lem-3.1}, we have
\begin{equation} \label{3.4-1}
	\begin{aligned}
	\sup_{t\in[0,T]}\left[
	\left\|(\theta-1)(t)\right\|_{L^2}^2
	+\left\|u(t)\right\|^4_{L^4}\right]
	+\int_0^T\int_{\mathbb{R}}\left[\frac{\kappa\theta_x^2}{v}
	+\theta\frac{\mu u_x^2}{v}+\frac{\mu u^2u_x^2}{v}\right]
	\lesssim F(\theta),
	\end{aligned}
\end{equation}
where
\begin{equation}\label{F-theta}
  F(\theta)=\left(1+\|\theta\|_{L^{\infty}_{t,x}}^{\frac{2}{1+2{\ell_1}}}\right)\left(1+\|\theta\|_{L^{\infty}_{t,x}}^{\frac1{{\ell_2}}}\right).
\end{equation}
\end{lemma}
\begin{proof}
In the same manner as Lemma \ref{lem_the1}, we define a set $$\Omega_a(t):=\{x\in\mathbb{R}:\ \theta(t,x)>a\}.$$
We now turn to the case $\chi=0$, which is considerably more involved than the case $\chi=1$. The proof will conclude with a brief discussion of the necessary modifications for $\chi=1$.

Firstly, multiplying \eqref{nsp}$_3$ and \eqref{nsp}$_2$ by $(\theta-2)_+:=\max\{\theta-2,0\}$ and $2u(\theta-2)_+$, respectively, adding the resultant together and integrating them over $[0,t]\times\mathbb{R}$, we derive that
\begin{align} \notag
	&\int_{\mathbb{R}}\left[\frac{c_v}2(\theta-2)_+^2+u^2(\theta-2)_+\right]\mathrm{d}x
	+\int_0^t\int_{\Omega_2(\tau)}\left[\frac{\kappa\theta_x^2}{v}
	+\frac{\mu u_x^2}{v}(\theta-2)_+\right] \\
	=&\int_{\mathbb{R}} \left[\frac{c_v}2(\theta_0-2)_+^2+u_0^2(\theta_0-2)_+\right]\mathrm{d}x
	+\underbrace{\int_0^t\int_{\mathbb{R}}\frac{R\theta}{v}(\theta-2)_+u_x}_{\mathcal{J}_1}+\underbrace{2\int_0^t\int_{\Omega_2(\tau)}\frac{R\theta}{v}u\theta_x}_{\mathcal{J}_2}\nonumber\\
    &-\underbrace{2\int_0^t\int_{\Omega_2(\tau)}\frac{\mu uu_x}{v}\theta_x}_{\mathcal{J}_3}
    +\underbrace{\frac1{c_v}\int_0^t\int_{\Omega_2(\tau)}u^2\left[\mu\frac{u_x^2}v-\frac{R\theta}vu_x\right]}_{\mathcal{J}_4}\nonumber\\
    &+\underbrace{\int_0^t\int_{\Omega_2(\tau)}\frac{1}{c_v}u^2\left(\frac{\kappa\theta_x}v\right)_x}_{\mathcal{J}_5}
    +\underbrace{2\int_0^t\int_{\mathbb{R}}\frac{\Phi_x}{v}u(\theta-2)_+}_{\mathcal{J}_6}\label{3.4-2}.
\end{align}
Different from the estimates in Lemma \ref{lem_the1}, we shall treat $\mathcal{J}_{p} (p=1\thicksim6)$ more careful. The following facts derived from \eqref{h} will be used frequently,
\begin{align}
	&\| h(v)^{-1}v^{-1}\|_{L^{\infty}_{t,x}}+\| h(v)^{-1}v\|_{L^{\infty}_{t,x}} \lesssim 1,\label{3.4-3}\\
	&\| \mu^{-1}v^{-1}\|_{L^{\infty}_{t,x}}+\| \mu^{-1} v\|_{L^{\infty}_{t,x}} \lesssim 1.\label{3.4-4}
\end{align}

We just collect the estimates of terms $\mathcal{J}_i$ for $i=1\sim5$ without proof for brevity, since it can be found in \cite{MR3564590}.
\begin{align}\label{J1-J5-2}
  \sum_{i=1}^5|\mathcal{J}_i|\leq \eta\left(\int_0^t\int_{\mathbb{R}}\frac{\mu u_x^2}{v}(\theta-2)_++\int_0^t\int_{\Omega_2(\tau)}\frac{\kappa\theta_x^2}{v}\right)
  +C(\eta)\left(\int_0^t\sup_{x\in\mathbb{R}}\left(\theta-\tfrac{3}{2}\right)_+^2
  +\int_0^t\int_{\Omega_2(\tau)}\frac{\mu u^2u_x^2}{v}\right).
\end{align}

For the term $\mathcal{J}_6$, we apply integration by parts and Lemma \ref{lem-2.2} to obtain
\begin{equation}\label{3.4-10}
  \mathcal{J}_6=-\int_{\mathbb{R}}\left(\frac{\Phi_x}{v}\right)^2(\theta-2)_++\int_{\mathbb{R}}\left(\frac{\Phi_{0x}}{v_0}\right)^2(\theta_0-2)_++\int_0^t\int_{\Omega_2(\tau)}\left(\frac{\Phi_x}{v}\right)^2\theta_t.
\end{equation}
Inserting \eqref{nsp}$_3$ into the last term on the right-hand side of \eqref{3.4-10}, we have
\begin{align*}
  \int_0^t\int_{\Omega_2(\tau)}\left(\frac{\Phi_x}{v}\right)^2\theta_t
  =\frac{1}{c_v}\int_0^t\int_{\Omega_2(\tau)}\left(\frac{\Phi_x}{v}\right)^2\left[\frac{\mu u_x^2}{v}-\frac{R\theta}{v}u_x+\left(\frac{\kappa\theta_x}{v}\right)_x\right]:=\sum_{p=1}^3\mathcal{K}_{p}.
\end{align*}

For the term $\mathcal{K}_{1}$, it derives from Sobolev's inequality, Cauchy's inequality and \eqref{basic-1} that
\begin{align}\notag
  |\mathcal{K}_{1}|&\lesssim \int_0^t\left\|\frac{\Phi_x}{v}\right\|_{L^{\infty}}^2\int_{\Omega_2(\tau)}\frac{\mu u_x^2}{v}\nonumber\\
  &\lesssim \int_0^t\left\|\frac{\Phi_x}{v}\right\|_{L^2}\|1-v\|_{L^2}\int_{\Omega_2(\tau)}\frac{\mu u_x^2}{v}\nonumber\\
  &\lesssim \int_0^t\int_{\Omega_2(\tau)}\frac{\mu u_x^2}{v}+\int_0^t\|v\|_{L^{\infty}}\int_{\Omega_2(\tau)}\frac{\mu u_x^2}{v}\nonumber\\
  &\leq \eta\int_0^t\int_{\Omega_2(\tau)}\frac{\mu\theta u_x^2}{v}+C(\eta)(1+\|v\|_{L^{\infty}_{t,x}}^2),\label{3.4-11}
\end{align}
where we have used the following facts derived by mean value theorem and Corollary \ref{cor-1},
\begin{equation}\label{3.4-K1}
  \int_{v\leq 2}|1-v|^2\leq \int_{\bbR}\phi(v)\leq C,\quad \int_{v\geq 2}\left|\frac{1-v}{v}\right|^2\leq \int_{\bbR}\phi(v)\leq C.
\end{equation}

For the term $\mathcal{K}_{2}$, using Cauchy's inequality and \eqref{basic-1} imply,
\begin{align}\notag
  |\mathcal{K}_{2}|\lesssim& \int_0^t\int_{\Omega_2(\tau)}\left[\left(\frac{\Phi_x}{v}\right)^2\theta^2+\frac{u_x^2}{v^2}\right]\\
  \lesssim& \int_0^t\sup_{x\in\mathbb{R}}\left(\theta-\tfrac{3}{2}\right)_+^2\mathrm{d}\tau+\|\mu^{-1}v^{-1}\|_{L^{\infty}_{t,x}}\int_0^t\int_{\Omega_2(\tau)}\frac{\mu u_x^2}{v}\nonumber\\
  \leq& C\int_0^t\sup_{x\in\mathbb{R}}\left(\theta-\tfrac{3}{2}\right)_+^2\mathrm{d}\tau+\eta\int_0^t\int_{\Omega_2(\tau)}\frac{\mu\theta u_x^2}{v}+C(\eta)\label{3.4-12}.
\end{align}

For the term $\mathcal{K}_{3}$, since the uniform upper and lower bounds of specific volume are not derived, we shall divide the integrant into two parts as follows,
\begin{align}\label{3.4-13}
  \mathcal{K}_{3}\leq \eta\int_0^t\int_{\Omega_2(\tau)}\frac{\kappa \theta_x^2}{v}+C(\eta)\int_0^t\left\|\frac{\Phi_x}{v}\right\|_{L^{2}}\left(\|1-v\|_{L^2(v\leq2)}+\|1-v\|_{L^2(v\geq2)}\right)\int_{\Omega_2(\tau)}\frac{\kappa (1-v)^2}{v}
\end{align}
Direct calculations show that
\begin{align}
 &\,\int_0^t\left\|\frac{\Phi_x}{v}\right\|_{L^2}\|1-v\|_{L^2(v\leq2)}\int_{\Omega_2(\tau)}\frac{\kappa (1-v)^2}{v}\nonumber\\
  \lesssim&\,  \int_0^t\int_{\Omega_2(\tau)}\left|g(v)(1-v)\frac{\kappa (1-v)}{vg(v)}\right|\nonumber\\
  \lesssim&\, \left\|\frac{\kappa (1-v)}{vg(v)}\right\|_{L^{\infty}_{t,x}}\int_0^t\int_{\Omega_2(\tau)}g(v)(v-1)\nonumber\\
  \lesssim&\, (1+\|v^{-1}\|_{L^{\infty}_{t,x}})\int_0^t\int_{\Omega_2(\tau)}g(v)(v-1)\label{3.4-14}
\end{align}
where we have used the fact that
\begin{equation*}
	\left|\frac{\kappa (1-v)}{vg(v)}\right|\sim
	\begin{cases}
	1,  & \quad {v\rightarrow\infty,}\\
	\frac{1}{v},  &\quad { v\rightarrow0,}
	\end{cases}
\end{equation*}
which yields
\begin{equation}\label{3.4-14-fact}
  \left\|\frac{\kappa (1-v)}{vg(v)}\right\|_{L^{\infty}_{t,x}}\lesssim 1+\|v^{-1}\|_{L^{\infty}_{t,x}}.
\end{equation}
Similarly,
\begin{align}\notag
  \int_0^t\left\|\frac{\Phi_x}{v}\right\|_{L^2}\|1-v\|_{L^2(v\geq2)}\int_{\Omega_2(\tau)}\frac{\kappa (1-v)^2}{v}\lesssim~ &
  \int_0^t\left\|\frac{\Phi_x}{v}\right\|_{L^2}\left\|\frac{1-v}{v}\right\|_{L^2(v\geq2)}\|v\|_{L^{\infty}}\int_{\Omega_2(\tau)}\frac{\kappa (1-v)^2}{v}\\
  \lesssim~ &\|v\|_{L^{\infty}_{t,x}}(1+\|v^{-1}\|_{L^{\infty}_{t,x}})\int_0^t\int_{\Omega_2(\tau)}g(v)(v-1)\label{3.4-15}.
\end{align}

On the other hand, the following fact holds,
\begin{align} 	\label{3.4-16}
	\int_0^t\int_{\mathbb{R}}\frac{\kappa\theta_x^2}{v}+\int_0^t\int_{\mathbb{R}}\theta\frac{\mu u_x^2}{v}
	\lesssim&~1+\int_0^t\int_{\Omega_2(\tau)}\frac{\kappa \theta_x^2}{v}+\int_0^t\int_{\Omega_2(\tau)}\frac{\mu u_x^2}{v}(\theta-2)_+.
\end{align}
Substituting \eqref{J1-J5-2}--\eqref{3.4-16} into \eqref{3.4-2}, we derive that
\begin{align}\notag
  &\int_{\mathbb{R}}\left[(\theta-2)_+^2+u^2(\theta-2)_++\left(\frac{\Phi_x}{v}\right)^2(\theta-2)_+\right]\mathrm{d}x
	+\int_0^t\int_{\mathbb{R}}\left[\frac{\kappa\theta_x^2}{v}
	+\theta\frac{\mu u_x^2}{v}\right] \\
	\lesssim~& \|v\|_{L^{\infty}_{t,x}}^2
    +(1+\|v\|_{L^{\infty}_{t,x}})(1+\|v^{-1}\|_{L^{\infty}_{t,x}})\int_0^t\int_{\bbR}g(v)(v-1)\nonumber\\
	&+\int_0^t\sup_{x\in\mathbb{R}}\left(\theta-\tfrac{3}{2}\right)_+^2\mathrm{d}\tau
	+\int_0^t\int_{\Omega_2(\tau)}\frac{\mu u^2u_x^2}{v}.\label{3.4-18}
\end{align}

To control the last term on the righthand-side of \eqref{3.4-18}, we multiply \eqref{nsp}$_2$ by $u^3$ and integrate the resultant over $[0,T]\times\bbR$ to get
\begin{align}\notag
  \frac{1}{4}\int_{\mathbb{R}}u^4dx
	+3\int_0^t\int_{\mathbb{R}}\frac{\mu u^2u_x^2}{v}
	=~&\frac{1}{4}\int_{\mathbb{R}}u_0^4dx+\underbrace{3R\int_0^t\int_{\bbR}u^2u_x\frac{1-v}{v}}_{\mathcal{I}_1}
    +\underbrace{3R\int_0^t\int_{\Omega_2(\tau)}u^2u_x\frac{\theta-1}{v}}_{\mathcal{I}_2}\\
    &+\underbrace{3R\int_0^t\int_{\mathbb{R}\setminus\Omega_2(\tau)}u^2u_x\frac{\theta-1}{v}}_{\mathcal{I}_3}
    +\underbrace{\int_0^t\int_{\bbR}\frac{\Phi_x}{v}u^3}_{\mathcal{I}_4}\label{3.4-19}.
\end{align}
We only show the estimates on $\mathcal{I}_4$, since $\mathcal{I}_{p}(p=1,2,3)$ can be calculated in the same way as Lemma 2.4 in \cite{MR3564590} by the second author.
For completeness, we list the results as follows,
\begin{align}\label{3.4-20}
  \sum_{i=1}^3|\mathcal{I}_i|\leq\eta\int_0^t\int_{\mathbb{R}}\frac{\mu\theta u_x^2}{v}+\eta\int_0^t\int_{\mathbb{R}}\frac{\mu u^2u_x^2}{v}
	+C(\eta)+C(\eta)\int_0^t\sup_{x\in\mathbb{R}}\left(\theta-\tfrac{3}{2}\right)_+^2\mathrm{d}\tau.
\end{align}

For the term $\mathcal{I}_4$, we derive from Lemma \ref{lem-2.2} that
\begin{align}\label{3.4-21}
  \mathcal{I}_4=-\int_{\mathbb{R}}\frac{1}{2}\left(\frac{\Phi_x}{v}\right)^2u^2+\int_{\mathbb{R}}\frac{1}{2}\left(\frac{\Phi_{0x}}{v_0}\right)^2u_0^2
  +\int_0^t\int_{\mathbb{R}}\left(\frac{\Phi_x}{v}\right)^2uu_t.
\end{align}
Insert \eqref{nsp}$_2$ into the last term on the righthand-side of \eqref{3.4-21} to get
\begin{align*}
  \int_0^t\int_{\mathbb{R}}\left(\frac{\Phi_x}{v}\right)^2uu_t=\int_0^t\int_{\mathbb{R}}\left(\frac{\Phi_x}{v}\right)^2u\left[-\left(\frac{R\theta}{v}\right)_x+\left(\frac{\mu u_x}{v}\right)_x+\frac{\Phi_x}{v}\right]
  :=\sum_{p=1}^3\mathcal{L}_{p}.
\end{align*}

For the term $\mathcal{L}_{1}$, we rewrite it as follows,
\begin{align}\notag
  \mathcal{L}_{1}=~&R\int_0^t\int_{\mathbb{R}}\left(\left(\frac{\Phi_x}{v}\right)^2u\right)_x\frac{1-v}{v}+R\int_0^t\int_{\Omega_2(\tau)}\left(\left(\frac{\Phi_x}{v}\right)^2u\right)_x\frac{\theta-1}{v}\\
  &+R\int_0^t\int_{\mathbb{R}\setminus\Omega_2(\tau)}\left(\left(\frac{\Phi_x}{v}\right)^2u\right)_x\frac{\theta-1}{v}:=\sum_{p=1}^3\mathcal{L}_{1p}.\label{3.4-L1}
\end{align}

For the term $\mathcal{L}_{11}$, we shall divide it into two parts for $v\leq 2$ and $v\geq2$,
\begin{align}\notag
  |\mathcal{L}_{11}|\lesssim~& \underbrace{\int_0^t\left\|\frac{\Phi_x}{v}\right\|_{L^2}\|1-v\|_{L^2}\left(\int_{v\leq2}\left|u_x\frac{1-v}{v}\right|+\int_{v\geq2}\left|u_x\frac{1-v}{v}\right|\right)}_{\mathcal{M}_1}\\
  &+\underbrace{\int_0^t\left\|\frac{\Phi_x}{v}\right\|_{L^{\infty}}\|u\|_{L^{\infty}}\left(\int_{v\leq 2}\frac{(1-v)^2}{v}+\int_{v\geq 2}\frac{(1-v)^2}{v}\right)}_{\mathcal{M}_2}.\label{3.4-22}
\end{align}
By using Cauchy's inequality, \eqref{basic-1}, \eqref{3.4-K1} and \eqref{3.4-14-fact}, we have
\begin{align}\notag
  \mathcal{M}_1\lesssim& \int_0^t\left\|\frac{\Phi_x}{v}\right\|_{L^2}\left\|1-v\right\|_{L^2}\left(\left\|\frac{u_x}{v}\right\|_{L^2(v\leq2)}\left\|1-v\right\|_{L^2(v\leq 2)}
  +\left\|u_x\right\|_{L^2(v\geq2)}\left\|\frac{1-v}{v}\right\|_{L^2(v\geq2)}\right)\nonumber\\
  \lesssim~& \int_0^t\left(\|1-v\|_{L^2}^2+\|u_x\|_{L^2}^2+\left\|\frac{u_x}{v}\right\|_{L^2}^2\right)\nonumber\\
  \lesssim~& \|\mu^{-1}v\|_{L^{\infty}_{t,x}}\left\|\frac{\mu (1-v)}{vg(v)}\right\|_{L^{\infty}_{t,x}}\int_0^t\int_{\mathbb{R}}g(v)(v-1)+\int_0^t\left(\|u_x\|_{L^2}^2+\left\|\frac{u_x}{v}\right\|_{L^2}^2\right)\nonumber\\
  \lesssim~& (1+\|v^{-1}\|_{L^{\infty}_{t,x}})\int_0^t\int_{\mathbb{R}}g(v)(v-1)+\eta\int_0^t\int_{\mathbb{R}}\frac{\mu\theta u_x^2}{v}+C(\eta)\label{3.4-23}
\end{align}
and
\begin{align}\notag
  \mathcal{M}_2\lesssim~& \int_0^t\left\|\frac{\Phi_x}{v}\right\|_{L^2}^{\frac{1}{2}}\left\|1-v\right\|_{L^2}^{\frac12}\left\|u\right\|_{L^2}^{\frac{1}{2}}\left\|u_x\right\|_{L^2}^{\frac{1}{2}}\left(\left\|\frac{1-v}{v}\right\|_{L^2(v\leq2)}+\left\|1-v\right\|_{L^2(v\geq2)}\right)\\
  \lesssim~& (\|\mu^{-1}v^{-1}\|_{L^{\infty}_{t,x}}+\|\mu^{-1}v\|_{L^{\infty}_{t,x}})\left\|\frac{\mu (1-v)}{vg(v)}\right\|_{L^{\infty}_{t,x}}\int_0^t\int_{\mathbb{R}}g(v)(v-1)+\int_0^t\|u_x\|_{L^2}^2\nonumber\\
  \lesssim~& (1+\|v^{-1}\|_{L^{\infty}_{t,x}})\int_0^t\int_{\mathbb{R}}g(v)(v-1)+\eta\int_0^t\int_{\mathbb{R}}\frac{\mu\theta u_x^2}{v}+C(\eta)\label{3.4-24}.
\end{align}
Hence, combining \eqref{3.4-23} and \eqref{3.4-24} yields
\begin{equation}\label{3.4-25}
  |\mathcal{L}_{11}|\lesssim(1+\|v^{-1}\|_{L^{\infty}_{t,x}})\int_0^t\int_{\mathbb{R}}g(v)(v-1)+\eta\int_0^t\int_{\mathbb{R}}\frac{\mu\theta u_x^2}{v}+C(\eta).
\end{equation}

The rest terms $\mathcal{L}_{12}$ and $\mathcal{L}_{13}$ can be estimated in the same way as $\mathcal{L}_{11}$,
\begin{align}\notag
  |\mathcal{L}_{12}|\lesssim~& \left|\int_0^t\int_{\Omega_2(\tau)}\left(\frac{\Phi_x}{v}\right)^2u_x\frac{\theta-1}{v}\right|
  +\left|\int_0^t\int_{\Omega_2(\tau)}\frac{\Phi_x}{v}(1-v)u\frac{\theta-1}{v}\right|\\
  \lesssim~& \int_0^t\int_{\Omega_2(\tau)}\left\{\left(\frac{\Phi_x}{v}\right)^2\left[\left(\frac{u_x}{v}\right)^2+\left(\theta-1\right)^2+\left(\frac{1-v}{v}\right)^2\right]+u^2(\theta-1)^2\right\}\nonumber\\
  \lesssim~& \eta\int_0^t\int_{\mathbb{R}}\frac{\mu\theta u_x^2}{v}+C(\eta)(1+\|v\|_{L^{\infty}_{t,x}}^2)+\int_0^t\sup_{x\in\mathbb{R}}\left(\theta-\tfrac{3}{2}\right)_+^2\nonumber\\
  &+(1+\|v^{-1}\|_{L^{\infty}_{t,x}})\int_0^t\int_{\mathbb{R}}g(v)(v-1)
  \label{3.4-26}
\end{align}
and
\begin{align}\notag
  |\mathcal{L}_{13}|\lesssim~& \left|\int_0^t\int_{\mathbb{R}\setminus\Omega_2(\tau)}\left(\frac{\Phi_x}{v}\right)^2u_x\frac{\theta-1}{v}\right|
  +\left|\int_0^t\int_{\mathbb{R}\setminus\Omega_2(\tau)}\frac{\Phi_x}{v}(1-v)u\frac{\theta-1}{v}\right|\\
  \lesssim~& \int_0^t\left\|\frac{\Phi_x}{v}\right\|_{L^2}\|1-v\|_{L^2}\left\|\frac{u_x}{v}\right\|_{L^2}\|\theta-1\|_{L^2(\mathbb{R}\setminus\Omega_2(\tau))}\nonumber\\
  &+\int_0^t\left\|\frac{\Phi_x}{v}\right\|_{L^2}^{\frac{1}{2}}\|1-v\|_{L^2}^{\frac{1}{2}}\|u\|_{L^2}^{\frac{1}{2}}\|u_x\|_{L^2}^{\frac{1}{2}}\left\|\frac{1-v}{v}\right\|_{L^2}\|\theta-1\|_{L^2(\mathbb{R}\setminus\Omega_2(\tau))}\nonumber\\
  \lesssim~& (1+\|v^{-1}\|_{L^{\infty}_{t,x}})\int_0^t\int_{\mathbb{R}}g(v)(v-1)+\eta\int_0^t\int_{\mathbb{R}}\frac{\mu\theta u_x^2}{v}+C(\eta)
  \label{3.4-27}.
\end{align}

Inserting \eqref{3.4-25}--\eqref{3.4-27} into \eqref{3.4-L1}, we obtain
\begin{align}\notag
  |\mathcal{L}_{1}|\lesssim~& (1+\|v^{-1}\|_{L^{\infty}_{t,x}})\int_0^t\int_{\mathbb{R}}g(v)(v-1)+\eta\int_0^t\int_{\mathbb{R}}\frac{\mu\theta u_x^2}{v}\\
  &+C(\eta)(1+\|v\|_{L^{\infty}_{t,x}}^2)
  +\int_0^t\sup_{x\in\mathbb{R}}\left(\theta-\tfrac{3}{2}\right)_+^2d\tau.\label{3.4-28}
\end{align}

For the term $\mathcal{L}_{2}$, we treat it similarly as $\mathcal{L}_{1}$ in the following way,
\begin{align}\notag
  |\mathcal{L}_{2}|\lesssim~& \left|\int_0^t\int_{\mathbb{R}}\left(\frac{\Phi_x}{v}\right)^2\frac{\mu u_x^2}{v}\right|
  +\left|\int_0^t\int_{\mathbb{R}}\frac{\Phi_x}{v}(1-v)u\frac{\mu u_x}{v}\right|\\
  \lesssim~& \eta\int_0^t\int_{\mathbb{R}}\frac{\mu u^2u_x^2}{v}+C(\eta)\int_0^t\int_{\mathbb{R}}\left(\frac{\Phi_x}{v}\right)^2\frac{\mu(1-v)^2}{v}+\int_0^t\int_{\bbR}\left(\frac{\Phi_x}{v}\right)^2\frac{\mu u_x^2}{v}\nonumber\\
  \lesssim~&  \eta\int_0^t\int_{\mathbb{R}}\frac{\mu u^2u_x^2}{v}+\eta\int_0^t\int_{\mathbb{R}}\frac{\mu\theta u_x^2}{v}+C(\eta)(1+\|v\|_{L^{\infty}_{t,x}}^2)\nonumber\\
  &+C(\eta)(1+\|v\|_{L^{\infty}_{t,x}})(1+\|v^{-1}\|_{L^{\infty}_{t,x}})\int_0^t\int_{\mathbb{R}}g(v)(v-1)\label{3.4-29}.
\end{align}

The term $\mathcal{L}_{3}$ directly satisfies
\begin{align}\label{3.4-30}
  \mathcal{L}_{3}=-\frac{1}{4}\int_{\bbR}\left(\frac{\Phi_x}{v}\right)^4+\frac{1}{4}\int_{\bbR}\left(\frac{\Phi_{0x}}{v_0}\right)^4.
\end{align}

Substituting \eqref{3.4-20}, \eqref{3.4-21} and \eqref{3.4-28}--\eqref{3.4-30} into \eqref{3.4-19}, we obtain
\begin{align}\notag
  &\int_{\mathbb{R}}\left[u^4+\left(\frac{\Phi_x}{v}\right)^4\right]dx
	+\int_0^t\int_{\mathbb{R}}\frac{\mu u^2u_x^2}{v}\nonumber\\
	\lesssim~& \eta\int_0^t\int_{\mathbb{R}}\frac{\mu\theta u_x^2}{v}+C(\eta)(1+\|v\|_{L^{\infty}_{t,x}}^2)+\int_0^t\sup_{x\in\mathbb{R}}\left(\theta-\tfrac{3}{2}\right)_+^2\nonumber\\
  &+C(\eta)(1+\|v\|_{L^{\infty}_{t,x}})(1+\|v^{-1}\|_{L^{\infty}_{t,x}})\int_0^t\int_{\mathbb{R}}g(v)(v-1)\label{3.4-31}.
\end{align}
In light of the fundamental theorem of calculus, we infer
\begin{align} \label{3.4-32}
	\int_0^t\sup_{x\in\mathbb{R}}\left(\theta-\tfrac{3}{2}\right)_+^2d\tau
	\lesssim~&~\eta\int_0^t\int_{\mathbb{R}}\frac{\kappa\theta_x^2}v
	+C(\eta).
\end{align}

Combining \eqref{3.4-18}, \eqref{3.4-31}--\eqref{3.4-32} and \eqref{vx-v}, we choose a suitably small $\eta$ such that
\begin{align}\notag
  &\int_{\mathbb{R}}\left[(\theta-1)^2+u^4+\left(\frac{\Phi_x}{v}\right)^4\right]dx
	+\int_0^t\int_{\mathbb{R}}\left[\frac{\kappa\theta_x^2}{v}
	+\theta\frac{\mu u_x^2}{v}+\frac{\mu u^2u_x^2}{v}\right] \\
  \lesssim~& 1+\|v\|_{L^{\infty}_{t,x}}^2+G(v)\int_0^t\int_{\mathbb{R}}g(v)(v-1), \label{3.4-33}
\end{align}
where
\begin{equation}\label{G-v}
  G(v):=(1+\|v\|_{L^{\infty}_{t,x}})(1+\|v^{-1}\|_{L^{\infty}_{t,x}}).
\end{equation}
In order to control $G(v)\int_0^t\int_{\mathbb{R}}g(v)(v-1)$, we use Young's inequality to rewrite \eqref{3.1-3} as follows
\begin{equation}\label{3.1-34}
  \left|\int_0^t\int_{\mathbb{R}}\frac{\mu u_x^2}{v}+\int_0^t\int_{\mathbb{R}}\frac{\mu v_x\theta_x}{v^2}\right|\leq \frac{1}{2}\int_0^t\int_{\mathbb{R}}\frac{\mu\theta v_x^2}{v^3}
  +\frac{1}{10G(v)}\int_0^t\int_{\mathbb{R}}\left[\frac{\kappa\theta_x^2}{v}
	+\theta\frac{\mu u_x^2}{v}\right]+CG(v),
\end{equation}
which yields
\begin{equation}\label{3.1-35}
  \sup_{t\in[0,T]}\left\|\frac{\mu v_x}{v}(t)\right\|_{L^2}^2
  +\int_0^T\int_{\mathbb{R}}\left[\frac{\mu \theta v_x^2}{v^3}+g(v)(v-1)\right]dxds
  \leq \frac{1}{10G(v)}\int_0^t\int_{\mathbb{R}}\left[\frac{\kappa\theta_x^2}{v}
	+\theta\frac{\mu u_x^2}{v}\right]+CG(v).
\end{equation}
Inserting \eqref{3.1-35} into \eqref{3.4-33}, we get
\begin{align}
  \int_{\mathbb{R}}\left[(\theta-1)^2+u^4+\left(\frac{\Phi_x}{v}\right)^4\right]dx
	+\int_0^t\int_{\mathbb{R}}\left[\frac{\kappa\theta_x^2}{v}
	+\theta\frac{\mu u_x^2}{v}+\frac{\mu u^2u_x^2}{v}\right]
\lesssim G(v)^2,\quad {\rm for}\ \chi=0. \label{3.4-36}
\end{align}

For the case of $\chi=1$, we only need to adjust the proof of $\mathcal{J}_6$ and $\mathcal{I}_4$ defined in \eqref{3.4-2} and \eqref{3.4-19}, respectively.
It derives from Sobolev's inequality, Cauchy's inequality and Corollary \ref{cor-1} that
\begin{align*}
  \int_0^t\int_{\mathbb{R}}\frac{\Phi_x}{v}u(\theta-2)_+\lesssim~& \int_0^t\int_{\mathbb{R}}\frac{\Phi^2_x}{v^2}+\int_0^t\int_{\mathbb{R}}u^2(\theta-2)^2_+\\
  \lesssim~& \|v^{-1}\|_{L^{\infty}_{t,x}}\int_0^T\int_{\bbR}\frac{\mu(e^{\Phi},\theta)}{v}\Phi^2_x+\int_{0}^T\sup_{x\in\mathbb{R}}(\theta-\frac32)_+^2,
\end{align*}
and
\begin{align*}
  \int_0^t\int_{\bbR}\frac{\Phi_x}{v}u^3\lesssim~& \int_0^t\int_{\mathbb{R}}\frac{\Phi^2_x}{v^2}+\int_0^t\int_{\mathbb{R}}u^6\\
  \lesssim~& 1+\|v^{-1}\|_{L^{\infty}_{t,x}}\int_0^T\int_{\bbR}\frac{\mu(e^{\Phi},\theta)}{v}\Phi^2_x+\eta\int_0^t\theta\frac{\mu u_x^2}{v}.
\end{align*}
The rest of the proof is the same. Finally, we obtain
\begin{align}
  \int_{\mathbb{R}}\left[(\theta-1)^2+u^4\right]dx
	+\int_0^t\int_{\mathbb{R}}\left[\frac{\kappa\theta_x^2}{v}
	+\theta\frac{\mu u_x^2}{v}+\frac{\mu u^2u_x^2}{v}\right]\lesssim G(v)^2,\quad {\rm for}\ \chi=1. \label{3.4-37}
\end{align}
Combining \eqref{3.4-36}, \eqref{3.4-37} and \eqref{v-bound}, we finish the proof.
\end{proof}

The following lemma is a direct conclusion of Lemma \ref{lem-3.4}.
\begin{lemma}\label{lem-3.5}
Under the assumptions in Lemma \ref{lem-3.1}, we have
\begin{equation}\label{vx-v-new}
  \sup_{t\in[0,T]}\left\|\frac{\mu v_x}{v}(t)\right\|_{L^2}^2
  +\int_0^T\int_{\mathbb{R}}\left(\frac{\mu \theta v_x^2}{v^3}+\chi\frac{\mu(e^{\chi\Phi},\theta)}{v}\Phi^2_x+\left(g(v,\theta)-g(e^{\chi\Phi},\theta)\right)\left(ve^{-\chi\Phi}-1\right)\right)
  \lesssim F(\theta)^{\frac{1}{2}}.
\end{equation}
\end{lemma}

Next, we present the first order estimates of $u(t,x)$.
\begin{lemma}\label{lem-3.6}
Under the assumptions in Lemma \ref{lem-3.1}, we have
\begin{align}\label{ux}
    \|u_x\|_{L^2}^2+\int_0^T\int_{\mathbb{R}}\frac{\mu u_{xx}^2}{v}
    \lesssim~& F(\theta)^2+\|\theta\|_{L^{\infty}_{t,x}}F(\theta)^{\frac{1}{2}}.
\end{align}
\end{lemma}
\begin{proof}
Multiplying \eqref{nsp}$_2$ by $u_{xx}$ and integrating the resultant over $[0,T]\times\bbR$, we apply Lemma \ref{lem2.2}, \eqref{3.4-4},
\eqref{h} and Sobolev's inequality to get
\begin{align}\notag
  &\|u_x\|^2+\|1-v\|^2+\int_0^T\int_{\mathbb{R}}\frac{\mu u_{xx}^2}{v}\\
    \lesssim~& 1+
	\left\|\frac{1}{\mu\kappa}\right\|_{L^{\infty}_{t,x}}
	\int_0^t\int_{\mathbb{R}}\frac{\kappa\theta_x^2}{v}
	+\left\|\frac{\theta}{\mu^2}\right\|_{L^{\infty}_{t,x}}\int_0^t\int_{\mathbb{R}}\frac{\mu\theta v_x^2}{v^3}\nonumber\\
	&+\left\|\frac{1}{\mu v}\right\|_{L^{\infty}_{t,x}}
	\int_0^t
	\left\|u_x\right\|_{L^\infty}^2\left\|\frac{\mu v_x}{v}\right\|_{L^2}^2
	+\int_{0}^t\int_{\mathbb{R}}\frac{\mu_x^2 u_x^2}{\mu v}\nonumber\\
    \lesssim~& F(\theta)+\|\theta\|_{L^{\infty}_{t,x}}F(\theta)^{\frac{1}{2}}+F(\theta)^2+\eta\int_0^t\int_{\mathbb{R}}\frac{\mu u_{xx}^2}{v},\quad {\rm for}\ \chi=0.\label{3.6-1}
\end{align}
The only difference between the cases $\chi=0$ and $\chi=1$ is the following estimates:
\begin{align}
  \int_0^T\int_{\mathbb{R}}\frac{\Phi_x}{v}u_{xx}
  \lesssim~& \eta\int_0^T\int_{\mathbb{R}}\frac{\mu u_{xx}^2}{v}+\int_0^T\int_{\mathbb{R}}\frac{\Phi_x^2}{\mu v}\nonumber\\
\lesssim~&  \eta\int_0^T\int_{\mathbb{R}}\frac{\mu u_{xx}^2}{v}+F(\theta)^{\frac{3}{2}},\quad {\rm for}\ \chi=1.\label{3.6-2}
\end{align}
Hence, we finish the proof.
\end{proof}

Now, we turn to derive the upper bound of the temperature $\theta(t,x)$.
\begin{lemma}\label{lem-3.7}
Under the assumptions in Lemma \ref{lem-3.1}, there exists positive constants $B_{i}(i=7,8)$ depend only the initial data, $R$ and $c_v$ such that
\begin{gather}  \label{upper_th}
	\theta(t,x) \leq B_7,\\[1.5mm]
	\label{bound_v}
	B_8^{-1}\leq v(t,x)\leq B_8,
	\\  		\label{E_1order}
	\|\theta_x\|_{L^2}^2
	+\int_0^t\int_{\bbR}\frac{\kappa\theta_{xx}^2}{v}ds\leq F(\theta)^4+\|\theta\|_{L^{\infty}_{t,x}}^2F(\theta).
\end{gather}
\end{lemma}
\begin{proof}
Multiplying \eqref{nsp}$_3$ by $\theta_{xx}$ and integrating the resultant over $[0,T]\times\bbR$, we deduce
\begin{align} \notag
	&\frac{c_v}{2}\|\theta_{x}(t)\|_{L^2}^2-\frac{c_v}{2}\|\theta_{0x}\|_{L^2}^2
	+\int_0^t\int_{\mathbb{R}}\frac{\kappa\theta_{xx}^2}{v}
	\\  &\qquad \label{3.7-1}
	=\int_0^t\int_{\mathbb{R}}
	\left[\frac{\kappa v_x\theta_x}{v^2}\theta_{xx}
	+\frac{\theta}{v}u_x\theta_{xx}
	-\frac{\mu u_x^2}{v}\theta_{xx}
	-\frac{\kappa_x\theta_x}{v}\theta_{xx}\right].
\end{align}
The four terms on the righthand-side of \eqref{3.7-1} can be estimated in the same way as Lemma 2.7 in \cite{MR3564590} by the second author, we just list the results for our case,
\begin{align} \label{3.7-2}
	&\left|\int_0^t\int_{\mathbb{R}}\frac{\kappa v_x\theta_x}{v^2}\theta_{xx}\right|+\left|\int_0^t\int_{\mathbb{R}}\frac{\theta}{v}u_x\theta_{xx}\right|
    +\left|\int_0^t\int_{\mathbb{R}}\frac{\mu u_x^2}{v}\theta_{xx}\right|+\left|\int_0^t\int_{\mathbb{R}}\frac{\kappa_x\theta_x}{v}\theta_{xx}\right|\nonumber\\	
	\lesssim~& \eta\int_0^t\int_{\mathbb{R}}\frac{\kappa\theta_{xx}^2}{v}
	+C(\eta)\left(F(\theta)^2+\|\theta\|_{L^{\infty}_{t,x}}F(\theta)+F(\theta)^4+\|\theta\|^2_{L^{\infty}_{t,x}}F(\theta)\right),
\end{align}
Inserting \eqref{3.7-2} into \eqref{3.7-1}, we get \eqref{E_1order} by choosing a suitably small $\eta$.

To derive the upper bound of the temperature, we combine \eqref{3.4-1} and \eqref{E_1order} to get
\begin{equation*}
  \|\theta-1\|_{L^{\infty}_{t,x}}^2\leq C\|\theta-1\|_{L^2}\|\theta_x\|_{L^2}\lesssim F(\theta)^{\frac{5}{2}}+\|\theta\|_{L^{\infty}_{t,x}}F(\theta).
\end{equation*}
Recall that
\begin{equation*}
  F(\theta)=\left(1+\|\theta\|_{L^{\infty}_{t,x}}^{\frac{2}{1+2{\ell_1}}}\right)\left(1+\|\theta\|_{L^{\infty}_{t,x}}^{\frac1{{\ell_2}}}\right).
\end{equation*}
Hence, we shall choose suitable $\ell_1$ and $\ell_2$ such that
\begin{equation}\label{ell-assump-1}
  \frac{2}{1+2\ell_1}+\frac{1}{\ell_2}< \frac{4}{5},
\end{equation}
then \eqref{upper_th} can be derived by Cauchy's inequality. Moreover, \eqref{bound_v} follows by inserting \eqref{upper_th} into \eqref{v-bound}.
\end{proof}

Combining Corollary \ref{cor-1} and Lemma \ref{lem-3.1}-\ref{lem-3.7}, we have
\begin{lemma}\label{lem-3.8}
Under the assumptions in Lemma \ref{lem-3.1}, there exists a constant $B_{i}(i=9,10)$ depends only on the initial data such that
\begin{align}\label{3.8-1}
|(v-1,u,\theta-1)(t)\|_1^2+\left\|\Phi_x\right\|_{1}^2
	+\int_0^t\left[\left\|\sqrt{\theta}v_x(s)\right\|_{L^2}^2
	+\|(\theta_x,u_x)(s)\|_1^2\right]\mathrm{d}s\leq B_9,\quad {\rm for}\ \chi=0,\\
|(v-1,u,\theta-1)(t)\|_1^2+\left\|\Phi\right\|_{2}^2
	+\int_0^t\left[\left\|\sqrt{\theta}v_x(s)\right\|_{L^2}^2
	+\|(u_x,\theta_x,\Phi_x)(s)\|_1^2\right]\mathrm{d}s\leq B_{10},\quad {\rm for}\ \chi=1.\label{3.8-2}
\end{align}
\end{lemma}

According to Lemma \ref{lem-3.8} and maximum principle, we can get the lower bound of the temperature.
\begin{lemma}[c.f.\cite{MR3564590}] \label{L_lower}
Under the assumptions in Lemma \ref{lem-3.1}, there exists a positive constant $B_{11}$ depending only on the initial data
such that
\begin{equation} \label{lower_th}
	\inf_{\mathbb{R}}\theta(t,\cdot)\geq \frac{\inf_{\mathbb{R}}\theta(s,\cdot)}
	{B_{11}\inf_{\mathbb{R}}\theta(s,\cdot)(t-s)+1}
	\quad
	\textrm{for all }
	0\leq s\leq t\leq T.
\end{equation}
\end{lemma}

We omit the proof on the estimates of higher order derivatives of $(v(t,x),u(t,x),\theta(t,x))$ whose procedures are similar to the Naiver--Stokes equations which can be found in \cite{MR3564590} by the second author.
\begin{corollary}[c.f.\cite{MR3564590}] \label{cor1}
Under the assumptions in Lemma \ref{lem-3.1}, there exists $C(m_2)>0$, which depends only on
$m_2$ and the initial data such that for all $t\in[0,T]$,
\begin{align}\label{higher-order}
 \|(v-1,u ,\theta-1)(t)\|_3^2+\left\|\Phi_x\right\|_{3}^2
	+\int_0^T
	\left[\|v_{x}(s)\|_2^2+\|(u_{x},\theta_{x})(s)\|_3^2\right]\mathrm{d}s
	\leq C(m_2),\quad {\rm for}\ \chi=0,\\
\|(v-1,u ,\theta-1)(t)\|_3^2+\left\|\Phi\right\|_{4}^2
	+\int_0^T
	\left[\|v_{x}(s)\|_2^2+\|(u_{x},\theta_{x},\Phi_x)(s)\|_3^2\right]\mathrm{d}s
	\leq C(m_2),\quad {\rm for}\ \chi=1.\label{higher-order-1}
\end{align}
\end{corollary}

\subsection{proof of Theorem \ref{main-theorem-1}}
With Lemma \ref{lem-3.7}, Lemma \ref{L_lower} and Corollary \ref{cor1} in hand, we can deduce Theorem \ref{main-theorem-1} by using the continuation argument introduced in \cite{MR3564590,MR3861291,MR3794193} and we omit the details for brevity.

\section{Appendix}
In the appendix, we will show some Lemmas which are used to prove the Theorem \ref{main-theorem} without proof.
Notice that the Boltzmann relation $n_e=e^{-\Phi}$ has brought both benefits and drawbacks.
In the following Lemma \ref{lem-2.2}, there exists a useful relationship pointwisely between $u(t,x)$ and $\Phi_{xt}(t,x)$ when $\chi=0$,
which is invalid when introducing the Boltzmann relation.
Instead, the Boltzmann relation can induce a damping effect which provides a uniform bound of $|\Phi(t,x)|$ and a good estimate for $\Phi_x(t,x)$.

First, we need the following lemma which is crucial to derive Lemma \ref{lem_the1} when $\chi=0$.
\begin{lemma}[c.f. \cite{MR3032988}]\label{lem-2.2}
Under the conditions in Lemma \ref{lem-2.1}, when $\chi=0$, for $\forall t \in [0, T ]$, we have
\begin{equation}
  u(t,x)=-\left(\frac{\Phi_x(t,x)}{v(t,x)}\right)_t.
\end{equation}
\end{lemma}

Next, we show the estimate of $\int_{\mathbb{R}}\frac{\Phi_x^2}{v}(0,x)dx$ when $\chi=1$, whose proof can be found in \cite{MR4755045,Zhao-Zhao-2026}.
\begin{lemma}\label{lem2.2}
	Under the assumptions of Theorem \ref{main-theorem} with $\chi=1$, there exist positive constants $C$ and $C(\underline{V},\overline{V})$ such that the following estimates hold for $\forall t \in [0, T ]$,
	\begin{align}\label{phi}
		\int_\mathbb{R}\frac{\Phi_x^2}{v}dx +	\int_\mathbb{R}v(1-\Phi e^{-\Phi}-e^{-\Phi})dx
		\leq C\left(\int_\mathbb{R}\frac{(v-1)^2}{v}dx+\left(\int_\mathbb{R}\frac{(v-1)^2}{v}dx\right)^3\right).
		\end{align}
	Consequently, one can deduce that
	\begin{align}\label{phi0}
		\left\|\Phi_{0x}\right\|^2+\left\|\sqrt{1-\Phi_0 e^{-\Phi_0}-e^{-\Phi_0}}\right\|^2\leq C(\underline{V},\overline{V})\left(\|v_0-1\|^2+\|v_0-1\|^6\right).
	\end{align}			
\end{lemma}

The next lemma focuses on the uniform bound on the electric potential $\Phi(t,x)$ when $\chi=1$, whose proof can be found in \cite{MR4755045,Zhao-Zhao-2026}.
\begin{lemma}\label{lem2.3}
		Under the assumptions of Theorem \ref{main-theorem} with $\chi=1$, there exists a generic positive constant $\overline{\Phi}$ which depends only on the initial data
 such that for $\forall t \in [0, T ]$,
		\begin{align*}
		|\Phi(t,x)|\leq \overline{\Phi}.
		\end{align*}
\end{lemma}
Notice that only Corollary \ref{cor_0.5} is needed to prove Lemma \ref{lem2.3} which means $\overline{\Phi}$ depends only on $R$, $c_v$ and $L^2$ norm of the initial perturbations.\\

\noindent{\bf Data Availability}\, Data sharing not applicable to this article as no data sets were generated or analysed during
the current study.\\

\noindent{\bf Declarations}\,
Conflict of interest on behalf of all authors, the corresponding author states that there is no conflict of
interest.



 \end{document}